\documentclass[reqno,12pt]{article}
\usepackage{amsmath}
\usepackage{amssymb}

\begin{document}

\renewcommand{\theequation}{\arabic{section}.\arabic{equation}}

\begin{center}
EXISTENCE FOR WEAKLY\ COERCIVE NONLINEAR DIFFUSION EQUATIONS VIA A
VARIATIONAL PRINCIPLE{\huge \ }

\bigskip

GABRIELA MARINOSCHI\footnote{%
Institute of Mathematical Statistics and Applied Mathematics of the Romanian
Academy, Calea 13 Septembrie 13, \ Bucharest, Romania. \noindent
gmarino@acad.ro, gabimarinoschi@yahoo.com
\par
{}}{\LARGE \ }

\bigskip

\bigskip
\end{center}

\noindent {\small Abstract. We are concerned with the study of the
well-posedness of a nonlinear diffusion equation with a monotonically
increasing multivalued time-dependent nonlinearity derived from a convex
continuous potential having a superlinear growth to infinity. The results in
this paper state that the solution of the nonlinear equation can be
retrieved as the null minimizer of an appropriate minimization problem for a
convex functional involving the potential and its conjugate. This approach,
inspired by the Brezis-Ekeland variational principle, provides new existence
results under minimal growth and coercivity conditions.}

\medskip


\medskip


\section{ Introduction}

\setcounter{equation}{0}

We are concerned with the study of the well-posedness of a nonlinear
diffusion equation with a monotonically increasing discontinuous
nonlinearity derived from a convex continuous potential, by using a dual
formulation of this equation as a minimization of an appropriate convex
functional. The idea of identifying the solutions of evolution equations as
the minima of certain functionals is due to Brezis and Ekeland and
originates in their papers published in 1976 (see \cite{brezis-ekeland-76}
and \cite{brezis-ekeland-76-2}). During the past decades this approach has
enjoyed much attention, as seen in the various literature and in some more
recently published monograph and papers (see e.g., \cite{auchmuty-88}, \cite%
{auchmuty-93}, \cite{ghoussoub}, \cite{ghoussoub-tzou-04}, \cite%
{ghoussoub-08}, \cite{visintin-08}, \cite{stefanelli-08}, \cite%
{stefanelli-09}, \cite{barbu-2011-jmma}, \cite{barbu-2012-jota}, \cite%
{gm-jota-var-1}). In \cite{gm-jota-var-1} two cases were considered, the
first for a continuous potential with a polynomial growth and the second for
a singular potential. The latter has provided the existence of the solution
to variational inequality which models a free boundary flow.

The challenging part in this duality principle is the proof of the
well-posedness of the evolution equation as a consequence of the existence
of a null minimizer in the associated minimization problem (that is a
solution which minimizes the functional to zero). A general receipt for
proving this implication does not exist, it rather depending on the good
choice of the functional and on the particularities of the potential of the
nonlinearity arising in the diffusion term. This way of approaching the
well-posedness of nonlinear diffusion equations by a dual formulation as a
minimization problem is extremely useful especially when a direct approach
by using the semigroup theory (see e.g., \cite{vb-springer-2010}, \cite%
{Crandall-Pazy-71}) or other classical variational results (see e.g., \cite%
{lions}) cannot be followed due either to the low regularity of the data or
to the weak coercivity of the potential.

In this work, the nonlinearity in the diffusion term is more general and it
has a time and space dependent potential assumed to have a weak coercivity
and no particular regularity with respect to time and space. The paper is
organized in two parts. At the beginning we investigate the case with the
potential and its conjugate depending on time and space. We prove that the
minimization problem has at least one solution, unique if the functional is
strictly convex. This seems to be a good candidate for the solution to the
nonlinear equation, reason for which it can be viewed as a generalized or
variational solution. If the admissible set is restricted by imposing a $%
L^{\infty }$-constraint on the state, then the generalized solution which
minimizes the functional to zero turns out to be quite the weak solution to
the nonlinear equation.

The second part concerns the case in which the potential does not depend on
space. The main result establishes that the null minimizer in the
minimization problem is the unique solution to the nonlinear equation,
provided that the potential exhibits a symmetry at large values of the
argument.

We would like to mention the benefit of such a duality approach, which
allows an elegant proof of the existence for a time dependent diffusion
equation, under general assumptions, by making possible its replacement by
the problem of minimizing a convex functional with a linear state equation.
We also stress that the existence results obtained in this way are not
covered and do not follow by the general existence theory of porous media
equations, as well as that of time dependent nonlinear infinite dimensional
Cauchy problems.

\section{Problem presentation}

We deal with the problem%
\begin{eqnarray}
\frac{\partial y}{\partial t}-\Delta \beta (t,x,y) &\ni &f\mbox{ \ \ \ \ \ \
\ \ \ \ \ \ \ \ in }Q:=(0,T)\times \Omega ,  \notag \\
-\frac{\partial \beta (t,x,y)}{\partial \nu } &=&\alpha \beta (t,x,y)\mbox{
\ \ on }\Sigma :=(0,T)\times \Gamma ,  \label{si-1} \\
y(0,x) &=&y_{0}\mbox{ \ \ \ \ \ \ \ \ \ \ \ \ \ in }\Omega ,  \notag
\end{eqnarray}%
where $\Omega $ is an open bounded subset of $\mathbb{R}^{N},$ $N\leq 3,$
with the boundary $\Gamma $ sufficiently smooth, $T$ is finite and $\beta $
has a potential $j$. The notation $\frac{\partial }{\partial \nu }$
represents the normal derivative and $\alpha $ is positive.

In this paper we assume that $j:Q\times \mathbb{R}\rightarrow (-\infty
,\infty ]$ and has the following properties:

\noindent $(h_{1})$ $\ (t,x)\rightarrow j(t,x,r)$ is measurable on $Q,$ for
all $r\in \mathbb{R},$

\noindent $(h_{2})$ $\ j(t,x,\cdot )$ is a proper, convex, continuous
function, a.e. $(t,x)\in Q,$
\begin{equation}
\partial j(t,x,r)=\beta (t,x,r)\mbox{ for all }r\in \mathbb{R},\mbox{ a.e. }%
(t,x)\in Q,  \label{si-beta-j}
\end{equation}%
\begin{eqnarray}
\frac{j(t,x,r)}{\left\vert r\right\vert } &\rightarrow &\infty ,\mbox{ \ as\
}\left\vert r\right\vert \rightarrow \infty ,\mbox{ uniformly for }(t,x)\in
Q,  \label{si-9-2} \\
\frac{j^{\ast }(t,x,\omega )}{\left\vert \omega \right\vert } &\rightarrow
&\infty ,\mbox{ \ as\ }\left\vert \omega \right\vert \rightarrow \infty ,%
\mbox{ uniformly for }(t,x)\in Q,  \label{si-9-2-00}
\end{eqnarray}%
\begin{equation}
j(\cdot ,\cdot ,0)\in L^{\infty }(Q),\mbox{ }j^{\ast }(\cdot ,\cdot ,0)\in
L^{\infty }(Q).  \label{si-9-2-0}
\end{equation}%
We define the conjugate $j^{\ast }:Q\times \mathbb{R}\rightarrow (-\infty
,\infty ]$ by
\begin{equation}
j^{\ast }(t,x,\omega )=\sup_{r\in \mathbb{R}}(\omega r-j(t,x,r)),\mbox{ a.e.
}(t,x)\in Q.  \label{si-4-0}
\end{equation}%
Then, the following two relations (Legendre-Fenchel) take place (see \cite%
{vb-springer-2010}, p. 6, see also \cite{fenchel-53}):
\begin{equation}
j(t,x,r)+j^{\ast }(t,x,\omega )\geq r\omega \mbox{ for all }r,\omega \in
\mathbb{R},\mbox{ a.e. }(t,x)\in Q,  \label{si-4-1}
\end{equation}%
\begin{equation}
j(t,x,r)+j^{\ast }(t,x,\omega )=r\omega ,\mbox{ iff }\omega \in \partial
j(t,x,r),\mbox{ for all }r,\omega \in \mathbb{R},\mbox{ a.e. }(t,x)\in Q.
\label{si-4-2}
\end{equation}

By (\ref{si-beta-j}) it follows that $\beta $ is a maximal monotone graph
(possibly multivalued) on $\mathbb{R}$, a.e. $(t,x)\in Q.$ Relations (\ref%
{si-9-2})-(\ref{si-9-2-00}) are equivalent with the the properties that $%
(\beta )^{-1}(t,x,\cdot )$ and $\beta (t,x,\cdot ),$ respectively, are
bounded on bounded subsets, uniformly a.e. $(t,x)\in Q$. This means that for
any $M>0$ there exists $Y_{M}$ and $W_{M},$ independent on $t$ and $x,$ such
that%
\begin{equation}
\sup \left\{ \left\vert r\right\vert ;\mbox{ }r\in \beta ^{-1}(t,x,\omega ),%
\mbox{ }\left\vert \omega \right\vert \leq M\right\} \leq W_{M},
\label{si-9-6}
\end{equation}%
\begin{equation}
\sup \left\{ \left\vert \omega \right\vert ;\mbox{ }\omega \in \beta (t,x,r),%
\mbox{ }\left\vert r\right\vert \leq M\right\} \leq Y_{M}.  \label{si-9-7}
\end{equation}%
In fact, when $j$ does not depend on $t$ and $x$, relations (\ref{si-9-2})-(%
\ref{si-9-2-00}) express\ that
\begin{equation*}
D(\partial j(r))=R(\partial j(r))=\mathbb{R}\mathbf{,}\mbox{ }D(\partial
j^{\ast }(r))=R(\partial j^{\ast }(r))=\mathbb{R}
\end{equation*}%
(see \cite{vb-springer-2010}, p. 9). We also recall that $\partial j^{\ast
}(t,x,\cdot )=(\partial j(t,x,\cdot ))^{-1}$ a.e. $(t,x)\in Q.$

We call \textit{weakly coercive }a nonlinear diffusion term with $j$ having
the properties (\ref{si-9-2})-(\ref{si-9-2-00}), and implicitly the
corresponding equation (\ref{si-1}).

We also recall that a proper, convex l.s.c. function is bounded below by an
affine function, hence
\begin{equation}
j(t,x,r)\geq k_{1}(t,x)r+k_{2}(t,x),\mbox{ }j^{\ast }(t,x,\omega )\geq
k_{3}(t,x)\omega +k_{4}(t,x)  \label{si-9-7-01}
\end{equation}%
for any $r,\omega \in \mathbb{R}$ and we assume that
\begin{equation}
k_{i}\in L^{\infty }(Q),\mbox{ }i=1,...4.  \label{si-9-7-02}
\end{equation}%
In fact (\ref{si-9-7-01}) follows if besides (\ref{si-9-2-0}) we assume that
there exist $\xi ,\eta \in L^{\infty }(Q)$ such that $\xi \in \partial
j(t,x,0),$ $\eta \in \partial j^{\ast }(t,x,0)$ a.e. $(t,x)\in Q.$

In this work we show that problem (\ref{si-1}) reduces to a certain
minimization problem $(P)$ for a convex lower semicontinuous functional
involving the functions $j$ and $j^{\ast }$. In Section 3, the existence of
at least a solution to $(P)$ is proved in Theorem 3.2, this being actually
the generalized solution associated to (\ref{si-1}). The uniqueness is
deduced directly from $(P)$ under the assumption of the strictly convexity
of $j.$ Moreover, when a state constraint $y\in \lbrack y_{m},y_{M}]$ is
included in the admissible set we show that the null minimization solution
is the unique weak solution to (\ref{si-1}) in Theorem 3.3.

In the case when $j$ does not depend on $x$ but on $t$ and has the same
behavior at $\left\vert r\right\vert $ large, i.e., it satisfies the
relation
\begin{equation}
j(t,-r)\leq \gamma _{1}j(t,r)+\gamma _{2},\mbox{ for any }r\in \mathbb{R},%
\mbox{ a.e. }t\in (0,T),  \label{si-9-5}
\end{equation}%
with $\gamma _{1}$ and $\gamma _{2}$ constants, we prove in Theorem 4.3 in
Section 4 that the solution to the minimization problem is the unique weak
solution to (\ref{si-1}) without assuming the previous additional state
constraint. This is based on Lemma 4.1 which plays an essential role in the
proof of this result. We mention that stochastic porous media equations of
the form (\ref{si-1}) were studied under a similar assumptions in \cite%
{barbu-daprato}, by a different method.

Theorem 4.3 is the main novelty of this work since it provides existence in (%
\ref{si-1}) for a time dependent weakly coercive $j$. With respect to the
treatment of the case which assumed a polynomial boundedness of $j$ (see
\cite{gm-jota-var-1}), the present one requires a sharp analysis in the $%
L^{1}$-space.

\subsection{Functional setting}

First, we introduce several linear operators related to problem (\ref{si-1}%
). Actually they represent the operator $-\Delta $ defined on various
spaces. The main operators which we use, $A_{0,\infty }$ and $A$ are defined
as follows:%
\begin{eqnarray}
A_{0,\infty }\psi &=&-\Delta \psi ,\mbox{ }A_{0,\infty }:D(A_{0,\infty
})=X\subset L^{\infty }(\Omega )\rightarrow L^{\infty }(\Omega ),
\label{200} \\
X &=&\left\{ \psi \in W^{2,\infty }(\Omega ),\mbox{ }\frac{\partial \psi }{%
\partial \nu }+\alpha \psi =0\mbox{ on }\Gamma \right\}  \notag
\end{eqnarray}%
and
\begin{eqnarray}
A &:&D(A)=L^{1}(\Omega )\subset X^{\prime }\rightarrow X^{\prime },  \notag
\\
\left\langle A\theta ,\psi \right\rangle _{X^{\prime },X} &=&\int_{\Omega
}\theta A_{0,\infty }\psi dx,\mbox{ }\forall \theta \in L^{1}(\Omega ),\mbox{
}\forall \psi \in X,  \label{A}
\end{eqnarray}%
where by $X^{\prime }$ we denote the dual of $X,$ with the pivot space $%
L^{2}(\Omega )$ ($X\subset L^{2}(\Omega )\subset X^{\prime }).$

We introduce the operator
\begin{eqnarray}
A_{1}\psi &=&-\Delta \psi ,\mbox{ }A_{1}:D(A_{1})\subset L^{1}(\Omega
)\rightarrow L^{1}(\Omega ),  \label{A0-L1} \\
D(A_{1}) &=&\left\{ \psi \in W^{1,1}(\Omega );\mbox{ }\Delta \psi \in
L^{1}(\Omega ),\mbox{ }\frac{\partial \psi }{\partial \nu }+\alpha \psi =0%
\mbox{ on }\Gamma \right\} ,  \notag
\end{eqnarray}%
which is $m$-accretive on $L^{1}(\Omega )$ (see \cite{brezis-strauss}). For
a later use we recall that
\begin{eqnarray}
A_{2}\psi &=&-\Delta \psi ,\mbox{ }A_{2}:X_{2}=D(A_{2})\subset L^{2}(\Omega
)\rightarrow L^{2}(\Omega ),  \label{A0} \\
X_{2} &=&\left\{ \psi \in W^{2,2}(\Omega );\mbox{ }\frac{\partial \psi }{%
\partial \nu }+\alpha \psi =0\mbox{ on }\Gamma \right\} ,  \notag
\end{eqnarray}%
is $m$-accretive on $L^{2}(\Omega )$ and $\widetilde{A_{2}},$ its extension
to $L^{2}(\Omega ),$ defined by%
\begin{eqnarray}
\widetilde{A_{2}} &:&L^{2}(\Omega )\subset X_{2}^{\prime }\rightarrow
X_{2}^{\prime },  \notag \\
\left\langle \widetilde{A_{2}}\theta ,\psi \right\rangle _{X_{2}^{\prime
},X_{2}} &=&\int_{\Omega }\theta A_{2}\psi dx,\mbox{ }\forall \theta \in
L^{2}(\Omega ),\mbox{ }\forall \psi \in X_{2},  \label{A2}
\end{eqnarray}%
is $m$-accretive on $X_{2}^{\prime }.$ Here, $X_{2}^{\prime }$ is the dual
of $X_{2}$ with $L^{2}(\Omega )$ as pivot space (see these last definitions
in \cite{gm-jota-var-1}).

Finally, let us consider the Hilbert space $V=H^{1}(\Omega )$ endowed with
the norm
\begin{equation*}
\left\Vert \phi \right\Vert _{V}=\left( \left\Vert \phi \right\Vert
^{2}+\alpha \left\Vert \phi \right\Vert _{L^{2}(\Gamma )}^{2}\right) ^{1/2},
\end{equation*}%
which is equivalent (for $\alpha >0)$ with the standard Hilbertian norm on $%
H^{1}(\Omega )$ (see \cite{necas-67}, p. 20)$.$ The dual of $V$ is denoted $%
V^{\prime }$ and the scalar product on $V^{\prime }$ is defined as%
\begin{equation}
(\theta ,\overline{\theta })_{V^{\prime }}=\left\langle \theta ,A_{V}^{-1}%
\overline{\theta }\right\rangle _{V^{\prime },V}  \label{si-0-0}
\end{equation}%
where $A_{V}:V\rightarrow V^{\prime }$ is given by%
\begin{equation}
\left\langle A_{V}\psi ,\phi \right\rangle _{V^{\prime },V}=\int_{\Omega
}\nabla \psi \cdot \nabla \phi dx+\int_{\Gamma }\alpha \psi \phi d\sigma ,%
\mbox{ for any }\phi \in V.  \label{si-0-1}
\end{equation}%
(In fact, $A_{V}$ is the extension of $A_{2}$ defined by (\ref{A0}) to $%
V^{\prime }).$

For the sake of simplicity, we shall omit sometimes to write the function
arguments in the integrands, writing $\int_{Q}gdxdt$ instead of $%
\int_{Q}g(t,x)dxdt,$ where $g:Q\rightarrow \mathbb{R}.$ In appropriate
places we indicate it as $g(t),$ to specify that $g:(0,T)\rightarrow Y,$
with $Y$ a Banach space.

\subsection{Statement of the problem}

In terms of the previously introduced operators we can write the abstract
Cauchy problem%
\begin{eqnarray}
\frac{dy}{dt}(t)+A\beta (t,x,y) &\ni &f(t),\mbox{ a.e. }t\in (0,T),
\label{si-7} \\
y(0) &=&y_{0}.  \notag
\end{eqnarray}

\medskip

\noindent \textbf{Definition 1.1.} Let $f\in L^{\infty }(Q)$ and $y_{0}\in
V^{\prime }.$ We call a\textit{\ weak solution} to (\ref{si-1}) a pair $%
(y,\eta ),$
\begin{equation*}
y\in L^{1}(Q)\cap W^{1,1}([0,T];X^{\prime }),\mbox{ }w\in L^{1}(Q),\mbox{ }%
w(t,x)\in \beta (t,x,y(t,x))\ \mbox{\ a.e. }(t,x)\in Q,
\end{equation*}%
which satisfies the equation
\begin{equation}
\int_{0}^{T}\left\langle \frac{dy}{dt}(t),\psi (t)\right\rangle _{X^{\prime
},X}dt+\int_{Q}w(t,x)(A_{0,\infty }\psi (t))(x)dxdt=\int_{0}^{T}\left\langle
f(t),\psi (t)\right\rangle _{X^{\prime },X}dt  \label{si-8-1-0}
\end{equation}%
for any $\psi \in L^{\infty }(0,T;X),$ and the initial condition $%
y(0)=y_{0}. $

\medskip

In literature, such a solution is called\ sometimes \textit{very weak} or
\textit{distributional} solution.

\medskip

We consider the minimization problem%
\begin{equation}
\underset{_{(y,w)\in U}}{\mbox{Minimize}}\mbox{ }J(y,w),  \tag{$P$}
\end{equation}%
where%
\begin{equation}
J(y,w)=\left\{
\begin{array}{l}
\int_{Q}\left( j(t,x,y(t,x))+j^{\ast }(t,x,w(t,x))\right) dxdt+\frac{1}{2}%
\left\Vert y(T)\right\Vert _{V^{\prime }}^{2} \\
\\
-\frac{1}{2}\left\Vert y_{0}\right\Vert _{V^{\prime
}}^{2}-\int_{Q}y(t,x)(A_{0,\infty }^{-1}f(t))(x)dxdt\mbox{\ \ if }(y,w)\in U,%
\mbox{ } \\
\\
+\infty ,\mbox{ \ \ \ \ \ \ \ \ \ \ \ \ \ \ \ \ \ \ \ \ \ \ \ \ \ \ \ \ \ \
\ \ \ \ \ \ \ \ \ \ \ \ \ \ \ \ \ otherwise,}%
\end{array}%
\right.  \label{J}
\end{equation}%
and%
\begin{multline*}
U=\{(y,w);\mbox{ }y\in L^{1}(Q)\cap W^{1,1}([0,T];X^{\prime }),\mbox{ }%
y(T)\in V^{\prime },\mbox{ }w\in L^{1}(Q), \\
j(\cdot ,\cdot ,y(\cdot ,\cdot ))\in L^{1}(Q),\mbox{ \ }j^{\ast }(\cdot
,\cdot ,w(\cdot ,\cdot ))\in L^{1}(Q), \\
(y,w)\mbox{ verifies (\ref{si-8-P}) below}\}
\end{multline*}%
\begin{eqnarray}
\frac{dy}{dt}(t)+Aw(t) &=&f(t)\mbox{ a.e. }t\in (0,T),  \label{si-8-P} \\
y(0) &=&y_{0}.  \notag
\end{eqnarray}%
Here, $\frac{dy}{dt}$ is taken in the sense of $X^{\prime }$-valued
distributions on $(0,T).$

We see that, by the existence theory of elliptic boundary value problems
(see \cite{ADN}), if $f(t)\in L^{\infty }(\Omega )$ then $A_{0,\infty
}^{-1}f(t)\in \bigcap\limits_{p\geq 2}W^{2,p}(\Omega )\subset L^{\infty
}(\Omega ),$ a.e. $t\in (0,T),$ so the last term in the expression of $J$
makes sense.

\section{Time and space dependent potential}

\setcounter{equation}{0}

In this section we consider that $j$ and $j^{\ast }$ depend on $t$ and $x$
as well, and assume $(h_{1})-(h_{2}),$ (\ref{si-beta-j})-(\ref{si-9-2-0}), (%
\ref{si-9-7-01})-(\ref{si-9-7-02}). We begin with an intermediate result.

\medskip

\noindent \textbf{Lemma 3.1. }\textit{The function }$J$\textit{\ is proper,
convex and lower semicontinuous on }$L^{1}(Q)\times L^{1}(Q).$

\medskip

\noindent \textbf{Proof.} It is obvious that\textbf{\ }$J$ is proper
(because $U\neq \varnothing )$ and convex. Let $\lambda >0.$ For the lower
semicontinuity we prove that the level set%
\begin{equation*}
E_{\lambda }=\{(y,w)\in L^{1}(Q)\times L^{1}(Q);\mbox{ }J(y,w)\leq \lambda \}
\end{equation*}%
is closed in $L^{1}(Q)\times L^{1}(Q).$ Let $(y_{n},w_{n})\in E_{\lambda }$
such that
\begin{equation}
y_{n}\rightarrow y\mbox{ strongly in }L^{1}(Q),\mbox{ \ }w_{n}\rightarrow w%
\mbox{ strongly in }L^{1}(Q),\mbox{ as }n\rightarrow \infty .
\label{si-198-1}
\end{equation}%
It follows that $(y_{n},w_{n})\in U$ is the solution to
\begin{eqnarray}
\frac{dy_{n}}{dt}(t)+Aw_{n}(t) &=&f(t),\mbox{ a.e. }t\in (0,T),
\label{si-closed-0} \\
y_{n}(0) &=&y_{0}  \notag
\end{eqnarray}%
and
\begin{eqnarray}
&&J(y_{n},w_{n})=\int_{Q}(j(t,x,y_{n}(t,x))+j^{\ast }(t,x,w_{n}(t,x)))dxdt
\label{si-198} \\
&&+\frac{1}{2}\left\{ \left\Vert y_{n}(T)\right\Vert _{V^{\prime
}}^{2}-\left\Vert y_{0}\right\Vert _{V^{\prime }}^{2}\right\}
-\int_{Q}y_{n}A_{0,\infty }^{-1}fdxdt\leq \lambda .\mbox{ }  \notag
\end{eqnarray}%
The convergences (\ref{si-198-1}) imply that
\begin{equation*}
\int_{0}^{T}\left\langle Aw_{n}(t),\psi (t)\right\rangle _{X^{\prime
},X}dt=\int_{Q}w_{n}A_{0,\infty }\psi dxdt\rightarrow \int_{Q}wA_{0,\infty
}\psi dxdt,\mbox{ as }n\rightarrow \infty ,
\end{equation*}%
for any $\psi \in L^{\infty }(0,T;X)$ and%
\begin{equation*}
\int_{Q}y_{n}A_{0,\infty }^{-1}fdxdt\rightarrow \int_{Q}yA_{0,\infty
}^{-1}fdxdt\mbox{, as }n\rightarrow \infty ,
\end{equation*}%
Therefore, by (\ref{si-closed-0}), we can write
\begin{equation}
\int_{0}^{T}\left\langle \frac{dy_{n}}{dt}(t),\psi (t)\right\rangle
_{X^{\prime },X}dt=-\int_{Q}w_{n}A_{0,\infty }\psi dxdt+\int_{Q}f\psi dxdt,
\label{si-198-2}
\end{equation}%
for any $\psi \in L^{\infty }(0,T;X),$ and we deduce that%
\begin{equation*}
\frac{dy_{n}}{dt}\rightarrow \frac{dy}{dt}\mbox{ weakly in }%
L^{1}(0,T;X^{\prime })\mbox{ as }n\rightarrow \infty ,
\end{equation*}%
meaning that $y_{n}$ is absolutely continuous on $[0,T]$ with values in $%
X^{\prime }.$

Again by (\ref{si-closed-0}) we have
\begin{equation}
y_{n}(t)=y_{0}+\int_{0}^{t}f(s)ds-\int_{0}^{t}Aw_{n}(s)ds\mbox{, for }t\in
\lbrack 0,T].  \label{si-200}
\end{equation}%
From here we get
\begin{equation}
\int_{\Omega }y_{n}(t)\phi dx=\left\langle y_{0},\phi \right\rangle
_{V^{\prime },V}+\int_{0}^{t}\int_{\Omega }f(s)\phi
dxds-\int_{0}^{t}\left\langle Aw_{n}(s),\phi \right\rangle _{X^{\prime },X}ds
\label{si-200-0}
\end{equation}%
for any $\phi \in X$ and $t\in \lbrack 0,T].$ Passing to the limit we obtain
\begin{equation*}
l(t)=\lim_{n\rightarrow \infty }\int_{\Omega }y_{n}(t)\phi dx=\left\langle
y_{0},\phi \right\rangle _{V^{\prime },V}+\int_{0}^{t}\int_{\Omega }f(s)\phi
dxds-\int_{0}^{t}\left\langle Aw(s),\phi \right\rangle _{X^{\prime },X}ds.
\end{equation*}%
We multiply this relation by $\varphi _{0}\in L^{\infty }(0,T)$ and
integrate over $(0,T),$ to obtain that%
\begin{eqnarray}
&&\int_{0}^{T}\varphi _{0}(t)l(t)dt  \label{si-201} \\
&=&\int_{0}^{T}\left( \left\langle y_{0},\phi \right\rangle _{V^{\prime
},V}+\int_{0}^{t}\int_{\Omega }f(s)\phi dxds-\int_{0}^{t}\left\langle
Aw(s),\phi \right\rangle _{X^{\prime },X}ds\right) \varphi _{0}(t)dt.  \notag
\end{eqnarray}%
We multiply (\ref{si-200}) by $\varphi _{0}(t)\phi (x)$ and integrate over $%
(0,T)\times \Omega .$ We have%
\begin{eqnarray}
\int_{Q}\varphi _{0}\phi y_{n}dxdt &=&\int_{0}^{T}\left( \left\langle
y_{0},\phi \right\rangle _{V^{\prime },V}+\int_{\Omega }\int_{0}^{t}f(s)\phi
dsdx\right) \varphi _{0}(t)dt  \label{si-202} \\
&&-\int_{0}^{T}\int_{0}^{t}\left\langle Aw_{n}(s),\phi \right\rangle
_{X^{\prime },X}\varphi _{0}(t)dsdt,  \notag
\end{eqnarray}%
whence we use the strong convergence $y_{n}\rightarrow y$ in $L^{1}(Q)$ to
get that
\begin{eqnarray}
\int_{Q}\varphi _{0}\phi ydxdt &=&\int_{0}^{T}\left( \left\langle y_{0},\phi
\right\rangle _{V^{\prime },V}+\int_{\Omega }\int_{0}^{t}f(s)\phi
dsdx\right) \varphi _{0}(t)dt  \label{si-203} \\
&&-\int_{0}^{T}\int_{0}^{t}\left\langle Aw(s),\phi \right\rangle _{X^{\prime
},X}\varphi _{0}(t)dsdt.  \notag
\end{eqnarray}%
Comparing (\ref{si-201}) and (\ref{si-203}) we deduce that
\begin{equation*}
\int_{0}^{T}\varphi _{0}(t)l(t)dt=\int_{Q}\varphi _{0}\phi ydxdt\mbox{ for
any }\varphi _{0}\in L^{\infty }(0,T),
\end{equation*}%
hence
\begin{equation*}
l(t)=\lim_{n\rightarrow \infty }\int_{\Omega }y_{n}(t)\phi dx=\int_{\Omega
}y(t)\phi dx\mbox{ for any }\phi \in X,\mbox{ }t\in \lbrack 0,T].
\end{equation*}%
Thus
\begin{equation}
y_{n}(t)\rightarrow y(t)\mbox{ weakly in }X^{\prime }\mbox{ as }n\rightarrow
\infty ,\mbox{ for any }t\in \lbrack 0,T]  \label{si-204}
\end{equation}%
and therefore
\begin{equation}
y_{n}(T)\rightarrow y(T)\mbox{, }y_{n}(0)\rightarrow y(0)=y_{0}\mbox{ weakly
in }X^{\prime }\mbox{, as }n\rightarrow \infty .  \label{si-206}
\end{equation}%
Letting $n\rightarrow \infty $ in (\ref{si-198-2}) we obtain
\begin{equation*}
\int_{0}^{T}\left\langle \frac{dy}{dt}(t),\psi (t)\right\rangle _{X^{\prime
},X}dt+\int_{0}^{T}\int_{\Omega }wA_{0,\infty }\psi
dxdt=\int_{0}^{T}\left\langle f(t),\psi (t)\right\rangle _{X^{\prime },X}dt,
\end{equation*}%
which proves that $(y,w)$ is the solution to (\ref{si-8-P}).

By (\ref{si-198}) and (\ref{si-9-7-01}) we can write that
\begin{eqnarray*}
&&\int_{Q}(k_{1}y_{n}+k_{2}+k_{3}w_{n}+k_{4})dxdt+\frac{1}{2}\left\Vert
y_{n}(T)\right\Vert _{V^{\prime }}^{2} \\
&\leq &\int_{Q}(j(t,x,y_{n}(t,x))+j^{\ast }(t,x,w_{n}(t,x)))dxdt+\frac{1}{2}%
\left\Vert y_{n}(T)\right\Vert _{V^{\prime }}^{2} \\
&\leq &\frac{1}{2}\left\Vert y_{0}\right\Vert _{V^{\prime }}^{2}+\left\Vert
y_{n}\right\Vert _{L^{1}(Q)}\left\Vert A_{0,\infty }^{-1}f\right\Vert
_{L^{\infty }(Q)}+\lambda \leq C\mbox{,}
\end{eqnarray*}%
whence, using (\ref{si-198-1}) we get
\begin{equation*}
\frac{1}{2}\left\Vert y_{n}(T)\right\Vert _{V^{\prime }}^{2}\leq
C+\max_{i=1,...4}\left\vert k_{i}\right\vert _{L^{\infty }(Q)}\left(
\left\Vert y\right\Vert _{L^{1}(Q)}+\left\Vert w\right\Vert
_{L^{1}(Q)}+2\right) =C_{1}
\end{equation*}%
with $C$ and $C_{1}$ constants and $\left\vert k_{i}\right\vert _{\infty
}=\left\Vert k_{i}\right\Vert _{L^{\infty }(Q)}.$

It follows that $y_{n}(T)\rightarrow \xi $ weakly in $V^{\prime }$ as $%
n\rightarrow \infty .$ As seen earlier, $y_{n}(T)\rightarrow y(T)$ weakly in
$X^{\prime },$ and by the uniqueness of the limit we get $\xi =y(T)\in
V^{\prime }.$

The function
\begin{equation*}
\varphi :L^{1}(Q)\rightarrow \mathbb{R},\mbox{ }\varphi
(z)=\int_{Q}j(t,x,z(t,x))dxdt
\end{equation*}%
is proper, convex and l.s.c. (see \cite{vb-springer-2010}, p. 56) and so by
Fatou's lemma (if $j$ would be nonnegative$)$ we get%
\begin{equation}
\varphi (y)\leq \liminf_{n\rightarrow \infty }\varphi
(y_{n})=\liminf_{n\rightarrow \infty }\int_{Q}j(t,x,y_{n}(t,x))dxdt<\infty .
\label{si-206-0}
\end{equation}%
Since $j$ is not generally nonnegative we use (\ref{si-9-7-01}) and apply
Fatou's lemma for
\begin{equation*}
\widetilde{j}(t,x,r)=j(t,x,r)-k_{1}(t,x)r-k_{2}(t,x)\geq 0.
\end{equation*}%
We get, by the strongly convergence $y_{n}\rightarrow y$ in $L^{1}(Q)$ and
the continuity of $j,$
\begin{eqnarray*}
&&\int_{Q}(j(t,x,y(t,x))-k_{1}y-k_{2})dxdt=\int_{Q}\liminf_{n\rightarrow
\infty }\widetilde{j}(t,x,y_{n}(t,x))dxdt \\
&\leq &\liminf_{n\rightarrow \infty }\int_{Q}\widetilde{j}%
(t,x,y_{n}(t,x))dxdt=\liminf_{n\rightarrow \infty
}\int_{Q}j(t,x,y_{n}(t,x))dxdt-\int_{Q}(k_{1}y+k_{2})dxdt,
\end{eqnarray*}%
and so (\ref{si-206-0}) holds.

Similarly we have that $\int_{Q}j^{\ast }(t,x,w(t,x))dxdt<\infty ,$ and so,
in particular, we have shown that $(y,w)\in U.$

Moreover, passing to the limit in (\ref{si-198}) as $n\rightarrow \infty $
we obtain by lower semicontinuity that
\begin{eqnarray*}
&&\int_{Q}(j(t,x,y(t,x))+j^{\ast }(t,x,w(t,x)))dxdt+\frac{1}{2}\left\Vert
y(T)\right\Vert _{V^{\prime }}^{2} \\
&&-\frac{1}{2}\left\Vert y_{0}\right\Vert _{V^{\prime
}}^{2}-\int_{Q}yA_{0,\infty }^{-1}fdxdt\leq \liminf_{n\rightarrow \infty
}J(y_{n},w_{n})\leq \lambda
\end{eqnarray*}%
which means that $(y,w)\in E_{\lambda }.$ This ends the proof. \ \hfill $%
\square $

\medskip

\medskip

\noindent \textbf{Theorem 3.2. }\textit{Problem }$(P)$\textit{\ has at least
a solution }$(y^{\ast },w^{\ast }).$ \textit{If} $j$ \textit{is strictly
convex the solution to }$(P)$\textit{\ is unique}.

\medskip

\noindent \textbf{Proof. }By (\ref{si-9-7-01}) we note that if $(y,w)\in U,$
then
\begin{eqnarray*}
J(y,w) &\geq &-\left\vert k_{1}\right\vert _{\infty }\left\Vert y\right\Vert
_{L^{1}(Q)}-\left\vert k_{2}\right\vert _{\infty }-\left\vert
k_{3}\right\vert _{\infty }\left\Vert w\right\Vert _{L^{1}(Q)}-\left\vert
k_{4}\right\vert _{\infty } \\
&&-\frac{1}{2}\left\Vert y_{0}\right\Vert _{V^{\prime }}^{2}-\left\Vert
y\right\Vert _{L^{1}(Q)}\left\Vert A_{0,\infty }^{-1}f\right\Vert
_{L^{\infty }(Q)}.
\end{eqnarray*}%
Let us set $d=\inf\limits_{(y,w)\in U}J(y,w).$ We assume first that $%
d>-\infty $ and we shall show later that this is indeed the only case.

Let us consider a minimizing sequence $(y_{n},w_{n})\in U$, such that
\begin{equation}
d\leq J(y_{n},w_{n})\leq d+\frac{1}{n},  \label{si-99}
\end{equation}%
where the pair $(y_{n},w_{n})$ satisfies (\ref{si-closed-0}).

By (\ref{si-9-2})-(\ref{si-9-2-00}), for any $M>0,$ there exist $C_{M}$ and $%
D_{M}$ such that $j(t,x,r)>M\left\vert r\right\vert $ as $\left\vert
r\right\vert >C_{M}$ and $j^{\ast }(t,x,\omega )>M\left\vert \omega
\right\vert $ as $\left\vert \omega \right\vert >D_{M}.$ Then, by (\ref%
{si-99}) we write%
\begin{eqnarray*}
&&\int_{\{(t,x);\left\vert y_{n}(t,x)\right\vert \leq
C_{M}\}}j(t,x,y_{n}(t,x))dxdt+M\int_{\{(t,x);\left\vert
y_{n}(t,x)\right\vert >C_{M}\}}\left\vert y_{n}\right\vert dxdt \\
&&+\int_{\{(t,x);\left\vert w_{n}(t,x)\right\vert \leq D_{M}\}}j^{\ast
}(t,x,w_{n}(t,x))dxdt+M\int_{\{(t,x);\left\vert w_{n}(t,x)\right\vert
>D_{M}\}}\left\vert w_{n}\right\vert dxdt \\
&&+\frac{1}{2}\left\Vert y_{n}(T)\right\Vert _{V^{\prime }}^{2}-\frac{1}{2}%
\left\Vert y_{0}\right\Vert _{V^{\prime }}^{2}\leq d+\frac{1}{n} \\
&\leq &\left\Vert A_{0,\infty }^{-1}f\right\Vert _{L^{\infty }(Q)}\left(
\int_{\{(t,x);\left\vert y_{n}(t,x\right\vert )\leq C_{M}\}}\left\vert
y_{n}\right\vert dxdt+\int_{\{(t,x);\left\vert y_{n}(t,x)\right\vert
>C_{M}\}}\left\vert y_{n}\right\vert dxdt\right) .
\end{eqnarray*}%
Denoting $\left\Vert A_{0,\infty }^{-1}f\right\Vert _{L^{\infty
}(Q)}=f_{\infty },$ and taking $M$ large enough such that $M>f_{\infty }$ it
follows that
\begin{eqnarray*}
&&(M-f_{\infty })\int_{\{(t,x);\left\vert y_{n}(t,x)\right\vert
>C_{M}\}}\left\vert y_{n}\right\vert dxdt+M\int_{\{(t,x);\left\vert
w_{n}(t,x)\right\vert >D_{M}\}}\left\vert w_{n}\right\vert dxdt \\
&&\mbox{ \ }+\frac{1}{2}\left\Vert y_{n}(T)\right\Vert _{V^{\prime }}^{2} \\
&\leq &\frac{1}{2}\left\Vert y_{0}\right\Vert _{V^{\prime }}^{2}+f_{\infty
}C_{M}\mbox{meas}(Q)+\int_{\{(t,x);\left\vert y_{n}(t,x)\right\vert \leq
C_{M}\}}\left\vert j(t,x,y_{n}(t,x))\right\vert dxdt \\
&&+\int_{\{(t,x);\left\vert w_{n}(t,x)\right\vert \leq D_{M}\}}\left\vert
j^{\ast }(t,x,w_{n}(t,x))\right\vert dxdt+d+1 \\
&\leq &\frac{1}{2}\left\Vert y_{0}\right\Vert _{V^{\prime }}^{2}+f_{\infty
}C_{M}\mbox{meas}(Q)+\int_{\{(t,x);\left\vert y_{n}(t,x)\right\vert \leq
C_{M}\}}\left\vert \widetilde{j}(t,x,y_{n}(t,x))\right\vert dxdt \\
&&+\int_{\{(t,x);\left\vert w_{n}(t,x)\right\vert \leq D_{M}\}}\left\vert
\widetilde{j}^{\ast }(t,x,w_{n}(t,x))\right\vert dxdt+d+1 \\
&&+\int_{\{(t,x);\left\vert y_{n}(t,x)\right\vert \leq C_{M}\}}\left\vert
k_{1}y_{n}+k_{2}\right\vert dxdt+\int_{\{(t,x);\left\vert
w_{n}(t,x)\right\vert \leq D_{M}\}}\left\vert k_{3}w_{n}+k_{4}\right\vert
dxdt,
\end{eqnarray*}%
where $\widetilde{j}(t,x,r)=j(t,x,r)-k_{1}r-k_{2},$ $\widetilde{j}^{\ast
}(t,x,\omega )=j^{\ast }(t,x,\omega )-k_{3}\omega -k_{4}.$

Recalling (\ref{si-9-7}) and (\ref{si-9-6}),
\begin{equation*}
j(t,x,y_{n}(t,x))\leq \left\vert j(t,x,0)\right\vert +\left\vert \eta
_{n}(t,x)\right\vert \left\vert y_{n}(t,x)\right\vert \leq Y_{M}^{1}\mbox{
on }\{(t,x);\left\vert y_{n}(t,x)\right\vert \leq C_{M}\},
\end{equation*}%
\begin{equation*}
j^{\ast }(t,x,w_{n}(t,x))\leq \left\vert j^{\ast }(t,x,0)\right\vert
+\left\vert \varpi _{n}(t,x)\right\vert \left\vert y(t,x)\right\vert \leq
W_{M}^{1}\mbox{ on }\{(t,x);\left\vert w_{n}(t,x)\right\vert \leq D_{M}\},
\end{equation*}%
where $\eta _{n}(t,x)\in \beta (t,x,y_{n})$ and $\varpi _{n}(t,x)\in (\beta
)^{-1}(t,x,w_{n})$ a.e. on $Q.$

Then
\begin{equation*}
0\leq \widetilde{j}(t,x,y_{n}(t,x))\leq Y_{M}^{1}+\left\vert
k_{1}\right\vert _{\infty }C_{M}+\left\vert k_{2}\right\vert _{\infty }\mbox{
on }\{(t,x);\left\vert y_{n}(t,x)\right\vert \leq C_{M}\},
\end{equation*}%
\begin{equation*}
0\leq \widetilde{j}^{\ast }(t,x,w_{n}(t,x))\leq W_{M}^{1}+\left\vert
k_{3}\right\vert _{\infty }D_{M}+\left\vert k_{4}\right\vert _{\infty }\mbox{
on }\{(t,x);\left\vert w_{n}(t,x)\right\vert \leq D_{M}\}
\end{equation*}%
and we deduce that%
\begin{eqnarray}
&&(M-f_{\infty })\int_{\{(t,x);\left\vert y_{n}(t,x)\right\vert
>C_{M}\}}\left\vert y_{n}\right\vert dxdt+M\int_{\{(t,x);\left\vert
w_{n}(t,x)\right\vert >D_{M}\}}\left\vert w_{n}\right\vert dxdt
\label{bound-inf} \\
&&+\frac{1}{2}\left\Vert y_{n}(T)\right\Vert _{V^{\prime }}^{2}\leq C+d.
\notag
\end{eqnarray}%
Consequently, this yields
\begin{equation}
\left\Vert y_{n}\right\Vert _{L^{1}(Q)}\leq C,\mbox{ \ }\left\Vert
w_{n}\right\Vert _{L^{1}(Q)}\leq C,\mbox{ \ }\left\Vert y_{n}(T)\right\Vert
_{V^{\prime }}\leq C.  \label{bound}
\end{equation}%
(By $C$ and $C_{i},$ $i=1,...4,$ we denote several constants independent on $%
n$).

From (\ref{si-99}) we get
\begin{equation}
I_{n}:=\int_{Q}j(t,x,y_{n}(t,x))dxdt+\int_{Q}j^{\ast
}(t,x,w_{n}(t,x))dxdt\leq C.  \label{si-198-0}
\end{equation}%
We continue by proving that separately each term is bounded, i.e.,
\begin{equation}
\int_{Q}j(t,x,y_{n}(t,x))dxdt\leq C_{1},\mbox{ }\int_{Q}j^{\ast
}(t,x,w_{n}(t,x))dxdt\leq C_{2}.  \label{si-100}
\end{equation}%
We write
\begin{eqnarray*}
&&I_{n}=\int_{\{(t,x);\left\vert y_{n}(t,x)\right\vert \leq
M\}}j(t,x,y_{n}(t,x))dxdt+\int_{\{(t,x);\left\vert y_{n}(t,x)\right\vert
>M\}}j(t,x,y_{n}(t,x))dxdt \\
&&+\int_{\{(t,x);\left\vert w_{n}(t,x)\right\vert \leq M\}}j^{\ast
}(t,x,w_{n}(t,x))dxdt+\int_{\{(t,x);\left\vert w_{n}(t,x)\right\vert
>M\}}j^{\ast }(t,x,w_{n}(t,x))dxdt \\
&\leq &C.
\end{eqnarray*}%
Therefore%
\begin{eqnarray*}
&&\int_{\{(t,x);\left\vert y_{n}(t,x)\right\vert
>M\}}j(t,x,y_{n}(t,x))dxdt+\int_{\{(t,x);\left\vert w_{n}(t,x)\right\vert
>M\}}j^{\ast }(t,x,w_{n}(t,x))dxdt \\
&\leq &C+Y_{M}^{1}\mbox{meas}(Q)+W_{M}^{1}\mbox{meas}(Q)=C_{3}.
\end{eqnarray*}%
Since $j(t,x,y_{n}(t,x))\geq k_{1}(t,x)y_{n}(t,x)+k_{2}(t,x)$ we deduce that%
\begin{equation*}
\int_{\{(t,x);\left\vert w_{n}(t,x)\right\vert >M\}}j^{\ast
}(t,x,w_{n}(t,x))dxdt\leq C_{4},
\end{equation*}%
whence
\begin{equation}
\int_{Q}j^{\ast }(t,x,w_{n}(t,x))dxdt\leq C_{1}.  \label{101-1}
\end{equation}%
Finally, (\ref{si-198-0}) yields
\begin{equation}
\int_{Q}j(t,x,y_{n}(t,x))dxdt\leq C_{2},  \label{101-2}
\end{equation}%
with $C_{1}$ and $C_{2}$ independent of $n.$

Next, we shall show that the sequences $(y_{n})_{n}$ and $(w_{n})_{n}$ are
weakly compact in $L^{1}(Q).$

To this end we have to show that the integrals $\int_{S}\left\vert
w_{n}\right\vert dxdt,$ with $S\subset Q,$ are equi-absolutely continuous,
meaning that for every $\varepsilon >0$ there exists $\delta $ such that $%
\int_{S}\left\vert w_{n}\right\vert dxdt<\varepsilon $ whenever meas$%
(S)<\delta .$ Let $M_{\varepsilon }>\frac{2C_{2}}{\varepsilon },$ where $%
C_{2}$ is the constant in (\ref{si-100}), and let $R_{M}$ be such that $%
\frac{j^{\ast }(t,x,w_{n})}{\left\vert w_{n}\right\vert }\geq M_{\varepsilon
}$ for $\left\vert r\right\vert >R_{M},$ by (\ref{si-9-2}). If $\delta <%
\frac{\varepsilon }{2R_{M}}$ then
\begin{eqnarray*}
&&\int_{S}\left\vert w_{n}\right\vert dxdt\leq \int_{\{(t,x);\left\vert
w_{n}(t,x)\right\vert >R_{M}\}}\left\vert w_{n}\right\vert
dxdt+\int_{\{(t,x);\left\vert w_{n}(t,x)\right\vert \leq R_{M}\}}\left\vert
w_{n}\right\vert dxdt \\
&\leq &M_{\varepsilon }^{-1}\int_{Q}j^{\ast
}(t,x,w_{n}(t,x))dxdy+R_{M}\delta <\varepsilon .
\end{eqnarray*}%
Hence, by the Dunford-Pettis theorem it follows that $(w_{n})_{n}$ is weakly
compact in $L^{1}(Q).$ In a similar way we proceed for showing the weakly
compactness of the sequence $(y_{n})_{n}.$ Thus,
\begin{equation*}
y_{n}\rightarrow y^{\ast }\mbox{ weakly in }L^{1}(Q),\mbox{ }%
w_{n}\rightarrow w^{\ast }\mbox{ weakly in }L^{1}(Q)\mbox{ as }n\rightarrow
\infty ,
\end{equation*}%
\begin{equation*}
Aw_{n}\rightarrow Aw^{\ast }\mbox{ weakly in }L^{1}(0,T;X^{\prime }),\mbox{
as }n\rightarrow \infty ,
\end{equation*}%
by (\ref{A}) which implies by (\ref{si-closed-0}) that
\begin{equation*}
\frac{dy_{n}}{dt}\rightarrow \frac{dy^{\ast }}{dt}\mbox{ weakly in }%
L^{1}(0,T;X^{\prime })\mbox{ as }n\rightarrow \infty .
\end{equation*}%
Passing to the limit in
\begin{equation*}
\int_{0}^{T}\left\langle \frac{dy_{n}}{dt}(t),\psi (t)\right\rangle
_{X^{\prime },X}dt+\int_{Q}w_{n}A_{0,\infty }\psi
dxdt=\int_{0}^{T}\left\langle f(t),\psi (t)\right\rangle _{X^{\prime },X}dt
\end{equation*}%
for any $\psi \in L^{\infty }(0,T;X)$ we get that $(y^{\ast },w^{\ast })$
verifies (\ref{si-8-1-0}), or equivalently (\ref{si-8-P}), i.e.,
\begin{equation*}
\int_{0}^{T}\left\langle \frac{dy^{\ast }}{dt}(t),\psi (t)\right\rangle
_{X^{\prime },X}dt+\int_{\Omega }w^{\ast }A_{0,\infty }\psi
dxdt=\int_{0}^{T}\left\langle f(t),\psi (t)\right\rangle _{X^{\prime },X}dt.
\end{equation*}%
Next we show that
\begin{equation*}
y_{n}(T)\rightarrow y^{\ast }(T)\mbox{ and }y_{n}(0)\rightarrow y(0)=y_{0}%
\mbox{ weakly in }V^{\prime }\mbox{, as }n\rightarrow \infty ,
\end{equation*}%
in a similar way as in Lemma 3.1. In order to obtain (\ref{si-203}) we use
the weakly compactness of ($y_{n})_{n}$ in $L^{1}(Q).$

Finally, by passing to the limit in (\ref{si-99}), on the basis of the
weakly lower semicontinuity of the functional $J$ on $L^{1}(Q)\times
L^{1}(Q),$ we obtain that
\begin{equation*}
J(y^{\ast },w^{\ast })=d.
\end{equation*}%
Hence, we have got that $y^{\ast }\in L^{1}(Q),$ $w^{\ast }\in L^{1}(Q),$ $%
y^{\ast }(T)\in V^{\prime }$ and $(y^{\ast },w^{\ast })$ satisfies (\ref%
{si-8-P}). By (\ref{si-100}) we get
\begin{equation*}
\int_{Q}j(t,x,y^{\ast }(t,x))dxdt<\infty ,\mbox{ }\int_{Q}j^{\ast
}(t,x,w^{\ast }(t,x))dxdt<\infty .
\end{equation*}%
With these relations we have ended the proof that $(y^{\ast },w^{\ast })$
belongs to $U$ and that it is is a solution to $(P).$

Let us show now that $d>-\infty .$ Indeed, otherwise, for every $K$ real
positive, there exists $n_{K},$ such that for every $n\geq n_{K}$ we have $%
J(y_{n},w_{n})<-K.$ Following the computations in the same way as before we
arrive at the inequality (\ref{bound-inf}) which reads now
\begin{eqnarray*}
&&(M-f_{\infty })\int_{\{(t,x);\left\vert y_{n}(t,x)\right\vert
>C_{M}\}}\left\vert y_{n}\right\vert dxdt+M\int_{\{(t,x);\left\vert
w_{n}(t,x)\right\vert >D_{M}\}}\left\vert w_{n}\right\vert dxdt \\
&&+\frac{1}{2}\left\Vert y_{n}(T)\right\Vert _{V^{\prime }}^{2}\leq C-K.
\end{eqnarray*}%
Since $C$ is a fixed constant, this implies $C-K<0,$ for $K$ large enough,
and this leads to a contradiction, as claimed.

The argument for the uniqueness proof is standard and it relies on the
assumption of the strict convexity of $j$ and on the obvious inequality
\begin{eqnarray}
&&J\left( \frac{y_{1}+y_{2}}{2},\frac{w_{1}+w_{2}}{2}\right)  \label{si-104}
\\
&=&\int_{Q}\left( j\left( t,x,\frac{y_{1}+y_{2}}{2}(t,x)\right) +j^{\ast
}\left( t,x,\frac{w_{1}+w_{2}}{2}(t,x)\right) \right) dxdt  \notag \\
&&+\frac{1}{2}\left\Vert \frac{y_{1}+y_{2}}{2}(T)\right\Vert _{V^{\prime
}}^{2}-\frac{1}{2}\left\Vert y_{0}\right\Vert _{V^{\prime }}^{2}-\int_{Q}%
\frac{y_{1}+y_{2}}{2}A_{0,\infty }^{-1}fdxdt  \notag \\
&\leq &\frac{1}{2}(J(y_{1},w_{1})+J(y_{2},w_{2}))-\frac{1}{2}\left\Vert
\frac{y_{1}-y_{2}}{2}(T)\right\Vert _{V^{\prime }}^{2},  \notag
\end{eqnarray}%
where $(y_{1},w_{1})$ and $(y_{2},w_{2})$ are two solutions to $(P).$\hfill $%
\square $

\medskip

We call the solution to the minimization problem $(P)$ a \textit{variational}
or \textit{generalized} solution to (\ref{si-1}).

One might suspect that if the minimum in $(P)$ is zero, then the null
minimizer is a weak solution to (\ref{si-1}). We shall prove this for a
slightly modified version of $(P),$ by including a boundedness constraint
for the state $y$ in the admissible set $U.$ More exactly we consider the
problem
\begin{equation}
\mbox{Minimize }\widetilde{J}(y,w)\mbox{ for all }(y,w)\in \widetilde{U}
\tag{$\widetilde{P}$}
\end{equation}%
where%
\begin{equation*}
\widetilde{J}(y,w)=\left\{
\begin{array}{l}
J(y,w),\mbox{ \ }(y,w)\in \widetilde{U}, \\
+\infty ,\mbox{ \ \ \ \ \ \ otherwise,}%
\end{array}%
\right.
\end{equation*}%
\begin{equation*}
\widetilde{U}=\{(y,w)\in U;\mbox{ }y(t,x)\in \lbrack y_{m},y_{M}]\mbox{ a.e.
}(t,x)\in Q\},
\end{equation*}%
with $y_{m},$ $y_{M}$ two constants. We assume that
\begin{equation*}
y_{0}\in L^{\infty }(\Omega ),\mbox{ }y_{0}\in \lbrack y_{m},y_{M}],\mbox{ \
}f\in L^{\infty }(Q)
\end{equation*}%
and remark that $\widetilde{U}$ is not empty (it contains e.g., $y_{0}$ with
$w_{0}=A_{0,\infty }^{-1}f(t))$ given by (\ref{si-8-P})).

If we set $y_{m}=0,$ then the previous boundedness property is in agreement
with the physical significance of $y,$ that of a fluid concentration in a
diffusion process, which is nonnegative.

Problem $(\widetilde{P})$ has at least a solution and the proof is the same
as in Theorem 3.2.

\medskip

\noindent \textbf{Theorem 3.3. }\textit{Let }$(y,w)\in \widetilde{U}$
\textit{be a null minimizer in }$(\widetilde{P}),$\textit{\ i.e., }%
\begin{equation*}
\min (\widetilde{P})=\widetilde{J}(y,w)=0.
\end{equation*}%
\textit{\ Let us assume in addition that }%
\begin{equation}
\frac{1}{2}\left\Vert y(T)\right\Vert _{V^{\prime }}^{2}-\frac{1}{2}%
\left\Vert y_{0}\right\Vert _{V^{\prime }}^{2}-\int_{Q}y(t,x)(A_{0,\infty
}^{-1}f(t))(x)dxdt=-\int_{Q}w(t,x)y(t,x)dxdt.  \label{si-298}
\end{equation}%
\textit{\ Then }%
\begin{equation*}
w(t,x)\in \beta (t,x,y(t,x)),\mbox{ }a.e.\mbox{ }(t,x)\in Q,
\end{equation*}%
\textit{and the pair }$(y,w)$\textit{\ is the unique weak solution to} (\ref%
{si-1}).

\medskip

\noindent \textbf{Proof. }Let $(y,w)$ be the null minimizer in $(\widetilde{P%
}).$ Then%
\begin{eqnarray*}
\widetilde{J}(y,w) &=&\int_{Q}(j(t,x,y(t,x))+j^{\ast }(t,x,w(t,x)))dxdt \\
&&+\frac{1}{2}\left\Vert y(T)\right\Vert _{V^{\prime }}^{2}-\frac{1}{2}%
\left\Vert y_{0}\right\Vert _{V^{\prime }}^{2}-\int_{Q}y(t,x)(A_{0,\infty
}^{-1}f(t))(x)dxdt=0.
\end{eqnarray*}%
By (\ref{si-298}) we have
\begin{equation}
\int_{Q}(j(t,x,y(t,x))+j^{\ast }(t,x,w(t,x))-y(t,x)w(t,x))dxdt=0.
\label{si-298-0}
\end{equation}%
This implies that $j(t,x,y(t,x))+j^{\ast }(t,x,w(t,x))-y(t,x)w(t,x)=0$ a.e. $%
(t,x)\in Q$ and so
\begin{equation*}
w(t,x)\in \beta (t,x,y(t,x))\mbox{ a.e. }(t,x)\in Q,
\end{equation*}%
as claimed. \hfill $\square $

\section{Time dependent potential}

\setcounter{equation}{0}

In this section we consider the case when $j$ and $j^{\ast }$ depend only on
$t$ and assume $(h_{1})-(h_{2}),$ (\ref{si-beta-j})-(\ref{si-9-2-0}), (\ref%
{si-9-7-01}) and (\ref{si-9-5}), where $k_{1},$ $k_{2}\in L^{\infty }(0,T).$

The main result of this section is that a solution to $(P)$ belongs to $%
L^{\infty }(0,T;V^{\prime })$ and minimizes $J$ to zero, being exactly the
unique weak solution to (\ref{si-1}).

To this end we need some intermediate results. The first is proved in the
next lemma and the second given in Theorem 4.2 recalls one of the main
results in \cite{gm-jota-var-1}.

\medskip

\noindent \textbf{Lemma 4.1. }\textit{Let} $(y,w)\in U$\textit{\ and }$y\in
L^{\infty }(0,T;V^{\prime }).$\textit{\ Then }$yw\in L^{1}(Q)$ \textit{and
we have the formula}
\begin{equation}
-\int_{Q}ywdxdt=\frac{1}{2}\left\Vert y(T)\right\Vert _{V^{\prime }}^{2}-%
\frac{1}{2}\left\Vert y_{0}\right\Vert _{V^{\prime
}}^{2}-\int_{Q}yA_{0,\infty }^{-1}fdxdt.  \label{60}
\end{equation}

\medskip

\noindent \textbf{Proof. }Let\textbf{\ }$(y,w)\in U.$ Then\ $y,w\in
L^{1}(Q), $ $j(\cdot ,\cdot ,y(\cdot ,\cdot ))\in L^{1}(Q),$ $j^{\ast
}(\cdot ,\cdot ,w(\cdot ,\cdot ))\in L^{1}(Q).$ By (\ref{si-9-5}) we have%
\begin{equation*}
j(t,-y(t,x))\leq \gamma _{1}j(t,y(t,x))+\gamma _{2}\mbox{ a.e. on }Q,
\end{equation*}%
which implies that $\int_{Q}j(t,x,-y(t,x))dxdt<\infty .$ Next, by the
relations%
\begin{eqnarray*}
j(t,x,y(t,x))+j^{\ast }(t,x,w(t,x)) &\geq &y(t,x)w(t,x), \\
j(t,x,-y(t,x))+j^{\ast }(t,x,w(t,x)) &\geq &-y(t,x)w(t,x)
\end{eqnarray*}%
it follows that
\begin{equation}
yw\in L^{1}(Q).  \label{61}
\end{equation}%
Because $(y,w)\in U$ it also satisfies (\ref{si-8-P}). Then $y\in
L^{1}(Q)\cap L^{\infty }(0,T;V^{\prime }),$ $w\in L^{1}(Q).$ We perform a
regularization by applying $(I+\varepsilon A_{\Delta })^{-1}$ to (\ref%
{si-8-P}), where $A_{\Delta }$ denotes here the realization of the operator $%
-\Delta $ on the spaces indicated in Section 2.1. We obtain%
\begin{eqnarray}
\frac{dy_{\varepsilon }}{dt}(t)+Aw_{\varepsilon }(t) &=&f_{\varepsilon }(t)%
\mbox{ a.e. }t\in (0,T),  \label{63} \\
y_{\varepsilon }(0) &=&(I+\varepsilon A_{V})^{-1}y_{0},  \notag
\end{eqnarray}%
where
\begin{eqnarray}
y_{\varepsilon }(t) &=&(I+\varepsilon A_{V})^{-1}y(t),\mbox{ a.e. }t\in (0,T)
\notag \\[1pt]
w_{\varepsilon }(t) &=&(I+\varepsilon A_{1})^{-1}w(t),\mbox{ a.e. }t\in (0,T)
\label{61-0} \\
f_{\varepsilon }(t) &=&(I+\varepsilon A_{0,\infty })^{-1}f(t),\mbox{ a.e. }%
t\in (0,T).  \notag
\end{eqnarray}%
According again to Brezis and Strauss (see \cite{brezis-strauss}), if $%
w(t)\in L^{1}(\Omega )$ then
\begin{equation}
w_{\varepsilon }(t)\in W^{1,q}(\Omega ),\mbox{ a.e. }t\in (0,T),\mbox{ with }%
1\leq q<\frac{N}{N-1}.  \label{64}
\end{equation}%
Since $\frac{N}{N-1}<N\leq 3$, we get by the Sobolev inequalities that
\begin{equation*}
W^{1,q}(\Omega )\subset L^{q^{\ast }}(\Omega ),\mbox{ }\frac{1}{q^{\ast }}=%
\frac{1}{q}-\frac{1}{N},
\end{equation*}%
with $\frac{N}{N-1}\leq q^{\ast }<\frac{N}{N-2}.$ It follows that
\begin{equation*}
w_{\varepsilon }\in L^{1}(0,T;L^{2}(\Omega )).
\end{equation*}%
Next,
\begin{equation*}
y_{\varepsilon }\in L^{1}(0,T;L^{2}(\Omega ))\cap L^{\infty }(0,T;V),
\end{equation*}%
by a similar argument as for $w_{\varepsilon },$ since $y\in L^{1}(Q)\cap
L^{\infty }(0,T:V^{\prime })$ and
\begin{equation*}
y_{\varepsilon }(t)=(I+\varepsilon A_{1})^{-1}y(t),\mbox{ a.e. }t\in (0,T),
\end{equation*}%
too. Finally,
\begin{equation*}
f_{\varepsilon }\in L^{\infty }(0,T;\bigcap\limits_{p\geq 2}W^{2,p}(\Omega
)),
\end{equation*}%
by the elliptic regularity.

Moreover, $A_{1}$ is $m$-accretive on $L^{1}(\Omega ),$ and it follows that
\begin{equation*}
w_{\varepsilon }(t)\rightarrow w(t)\mbox{ strongly in }L^{1}(\Omega )\mbox{
for any }t\in \lbrack 0,T]
\end{equation*}%
and
\begin{equation*}
\left\Vert w_{\varepsilon }(t)\right\Vert _{L^{1}(\Omega )}\leq \left\Vert
w(t)\right\Vert _{L^{1}(\Omega )}\mbox{ for any }t\in \lbrack 0,T],
\end{equation*}%
(see \cite{brezis-strauss}).

For a later use, we deduce by the Lebesgue dominated convergence theorem
that
\begin{equation}
w_{\varepsilon }\rightarrow w\mbox{ strongly in }L^{1}(Q),\mbox{ as }%
\varepsilon \rightarrow 0.  \label{63-0}
\end{equation}%
Similarly, we have that
\begin{equation}
y_{\varepsilon }\rightarrow y\mbox{ strongly in }L^{1}(Q),\mbox{ as }%
\varepsilon \rightarrow 0.  \label{63-1}
\end{equation}%
Finally,%
\begin{equation}
f_{\varepsilon }\rightarrow f\mbox{ weak* in }L^{\infty }(Q),\mbox{ and
strongly in }L^{p}(Q),\mbox{ }p\geq 2,\mbox{ as }\varepsilon \rightarrow 0.
\label{63-2}
\end{equation}%
By the first relation in (\ref{61-0}) we still have that
\begin{equation*}
(I+\varepsilon A_{V})^{-1}y(t)\rightarrow y(t)\mbox{ strongly in }V^{\prime }%
\mbox{ for any }t\in \lbrack 0,T].
\end{equation*}%
We also observe that
\begin{equation*}
\int_{0}^{T}\left\langle Aw_{\varepsilon }(t),\psi (t)\right\rangle
_{X^{\prime },X}dt=\int_{Q}w_{\varepsilon }A_{0,\infty }\psi dxdt\rightarrow
\int_{Q}wA_{0,\infty }\psi dxdt\mbox{ as }\varepsilon \rightarrow 0,
\end{equation*}%
for any $\psi \in L^{\infty }(0,T;X)$ and by (\ref{63})%
\begin{equation*}
\frac{dy_{\varepsilon }}{dt}\rightarrow \frac{dy}{dt}\mbox{ weakly in }%
L^{1}(0,T;X^{\prime })\mbox{ as }\varepsilon \rightarrow 0.
\end{equation*}%
Passing to the limit in (\ref{63}) tested for any $\psi \in L^{\infty
}(0,T;X),$%
\begin{equation*}
\int_{0}^{T}\left\langle \frac{dy_{\varepsilon }}{dt}(t),\psi
(t)\right\rangle _{X^{\prime },X}dt+\int_{Q}w_{\varepsilon }A_{0,\infty
}\psi dxdt=\int_{0}^{T}\left\langle f_{\varepsilon }(t),\psi
(t)\right\rangle _{X^{\prime },X}dt
\end{equation*}%
we check that $(y,w)$ indeed satisfies (\ref{si-8-P}).

Next, we assert that
\begin{equation}
\int_{Q}j(t,y_{\varepsilon }(t,x))dxdt\leq \int_{Q}j(t,y(t,x))dxdt.
\label{68}
\end{equation}

Indeed, let us introduce the Yosida approximation of $\beta ,$
\begin{equation}
\beta _{\lambda }(t,r)=\frac{1}{\lambda }\left( 1-(1+\lambda \beta (t,\cdot
)\right) ^{-1})r,\mbox{ a.e. }t,\mbox{ for all }r\in \mathbb{R}\mbox{ and }%
\lambda >0.  \label{69}
\end{equation}%
We have $\beta _{\lambda }(t,r)=\frac{\partial j_{\lambda }}{\partial r}%
(t,r),$ where $j_{\lambda }$ is the Moreau approximation of $j,$
\begin{equation}
j_{\lambda }(t,r)=\inf_{s\in \mathbb{R}}\left\{ \frac{\left\vert
r-s\right\vert ^{2}}{2\lambda }+j(t,s)\right\} ,\mbox{ a.e. }t,\mbox{ for
all }r\in \mathbb{R},  \label{70}
\end{equation}%
that can be still written as
\begin{equation}
j_{\lambda }(t,r)=\frac{1}{2\lambda }\left\vert (1+\lambda \beta (t,\cdot
))^{-1}r-r\right\vert ^{2}+j(t,(1+\lambda \beta (t,\cdot ))^{-1}r)).
\label{70-0}
\end{equation}%
The function $j_{\lambda }$ is convex, continuous and satisfies
\begin{eqnarray}
j_{\lambda }(t,r) &\leq &j(t,r)\mbox{ for all }r\in \mathbb{R},\mbox{ }%
\lambda >0,\mbox{ }  \label{71} \\
\lim_{\lambda \rightarrow 0}j_{\lambda }(t,r) &=&j(t,r),\mbox{ for all }r\in
\mathbb{R}.  \notag
\end{eqnarray}%
We have
\begin{eqnarray*}
&&\int_{Q}j_{\lambda }(t,y_{\varepsilon }(t,x))dxdt \\
&\leq &\int_{Q}j_{\lambda }(t,y(t,x))dxdt-\varepsilon
\int_{0}^{T}\left\langle A_{V}y_{\varepsilon }(t),\beta _{\lambda
}(t,y_{\varepsilon }(t))\right\rangle _{V^{\prime },V}dt.
\end{eqnarray*}%
Since for any $z\in V$ one has
\begin{equation*}
-\left\langle A_{V}z,\beta _{\lambda }(t,z)\right\rangle _{V^{\prime
},V}=-\int_{\Gamma }\alpha \beta _{\lambda }(t,z)zd\sigma -\int_{\Omega }%
\frac{\partial \beta _{\lambda }}{\partial z}(t,z)\left\vert \nabla
z\right\vert ^{2}dx\leq 0,
\end{equation*}%
we obtain
\begin{equation}
\int_{Q}j_{\lambda }(t,y_{\varepsilon }(t,x))dxdt\leq \int_{Q}j_{\lambda
}(t,y(t,x))dxdt.  \label{70-1}
\end{equation}%
Now, by (\ref{71})%
\begin{equation*}
j_{\lambda }(t,y(t,x))\leq j(t,y(t,x))
\end{equation*}%
and
\begin{equation}
\lim_{\lambda \rightarrow 0}j_{\lambda }(t,y(t,x))=j(t,y(t,x))\mbox{ a.e. on
}Q.  \label{70-2}
\end{equation}%
Let us assume for the moment that $j$ would be nonnegative. Then by the
Lebesgue dominated convergence theorem we get
\begin{equation}
\lim_{\lambda \rightarrow 0}\int_{Q}j_{\lambda
}(t,y(t,x))dxdt=\int_{Q}j(t,y(t,x))dxdt,\mbox{ for any }y\mbox{ fixed.}
\label{70-3}
\end{equation}%
Passing to the limit in (\ref{70-1}) as $\lambda \rightarrow 0,$ we obtain
that
\begin{equation}
\int_{Q}j(t,y_{\varepsilon }(t,x))dxdt\leq \int_{Q}j(t,y(t,x))dxdt,\mbox{
for all }\varepsilon >0.  \label{70-4}
\end{equation}%
Since in general $j$ is not necessarily nonnegative we consider the function
\begin{equation}
\widetilde{j}(t,r)=j(t,r)-k_{1}(t)r-k_{2}(t)  \label{70-6}
\end{equation}%
which is nonnegative by (\ref{si-9-7-01}). Hence, by (\ref{70-4}) we have
\begin{equation*}
\int_{Q}\widetilde{j}(t,y_{\varepsilon }(t,x))dxdt\leq \int_{Q}\widetilde{j}%
(t,y(t,x))dxdt+\int_{Q}k_{1}(t)(y(t,x)-y_{\varepsilon }(t,x))dxdt,
\end{equation*}%
hence
\begin{equation}
\int_{Q}\widetilde{j}(t,y_{\varepsilon }(t,x))dxdt\leq \int_{Q}\widetilde{j}%
(t,y(t,x))dxdt+\delta (\varepsilon ),  \label{70-7}
\end{equation}%
\begin{equation*}
\delta (\varepsilon )=\int_{Q}k_{1}(t)(y(t,x)-y_{\varepsilon }(t,x))dxdt\leq
\left\Vert k_{1}\right\Vert _{L^{\infty }(0,T)}\left\Vert y-y_{\varepsilon
}\right\Vert _{L^{1}(Q)}\rightarrow 0,\mbox{ as }\varepsilon \rightarrow 0,
\end{equation*}%
by (\ref{63-1}). Then, (\ref{70-7}) implies (\ref{68}) as claimed.

A similar relation to (\ref{68}) takes place for $j^{\ast }$,%
\begin{equation}
\int_{Q}j^{\ast }(t,w_{\varepsilon }(t,x))dxdt\leq \int_{Q}j(t,w(t,x))dxdt.
\label{70-5}
\end{equation}

This implies that $j(\cdot ,y_{\varepsilon }(\cdot ,\cdot ))\in L^{1}(Q)$, $%
j^{\ast }(\cdot ,w_{\varepsilon }(\cdot ,\cdot ))\in L^{1}(Q),$ for all $%
\varepsilon >0,$ and so, by the same argument as for $yw$ we deduce that
\begin{equation*}
y_{\varepsilon }w_{\varepsilon }\in L^{1}(Q).
\end{equation*}

We test (\ref{63}) by $A_{2}^{-1}y_{\varepsilon }(t)$ and integrate over $%
(0,T).$ Since $y_{\varepsilon }\in L^{1}(Q)\cap L^{\infty }(0,T;V),$ $%
w_{\varepsilon }\in L^{1}(0,T;L^{2}(\Omega ))$ we get $A_{2}^{-1}y_{%
\varepsilon }(t)\in X_{2},$ a.e. $t\in (0,T)$ and by (\ref{A2})
\begin{equation*}
\int_{0}^{T}\left\langle \widetilde{A_{2}}w_{\varepsilon
}(t),A_{2}^{-1}y_{\varepsilon }(t)\right\rangle _{X_{2}^{\prime
},X_{2}}dt=\int_{Q}y_{\varepsilon }w_{\varepsilon }dxdt.
\end{equation*}%
Then, by a few computations we deduce by (\ref{63}) that
\begin{equation}
-\int_{Q}y_{\varepsilon }w_{\varepsilon }dxdt=\frac{1}{2}\left\Vert
y_{\varepsilon }(T)\right\Vert _{V^{\prime }}^{2}-\frac{1}{2}\left\Vert
(I+\varepsilon A_{V})^{-1}y_{0}\right\Vert _{V^{\prime
}}^{2}-\int_{Q}y_{\varepsilon }A_{0,\infty }^{-1}f_{\varepsilon }dxdt.
\label{65}
\end{equation}%
Recalling that by (\ref{61-0}) we have that
\begin{equation*}
(I+\varepsilon A_{V})^{-1}y(t)\rightarrow y(t)\mbox{ strongly in }V^{\prime }%
\mbox{ for any }t\in \lbrack 0,T]
\end{equation*}%
and passing to the limit in (\ref{65}) as $\varepsilon \rightarrow 0$ we
obtain
\begin{equation}
\lim_{\varepsilon \rightarrow 0}\left( -\int_{Q}y_{\varepsilon
}w_{\varepsilon }dxdt\right) =\frac{1}{2}\left\Vert y(T)\right\Vert
_{V^{\prime }}^{2}-\frac{1}{2}\left\Vert y_{0}\right\Vert _{V^{\prime
}}^{2}-\int_{Q}yA_{0,\infty }^{-1}fdxdt.  \label{66}
\end{equation}%
Moreover, by the strongly convergence of $(y_{\varepsilon }$)$_{\varepsilon
} $ and $(w_{\varepsilon })_{\varepsilon },$ (\ref{63-1}) and (\ref{63-0})
we get
\begin{equation*}
y_{\varepsilon }\rightarrow y\mbox{ a.e. in }Q,\mbox{ }w_{\varepsilon
}\rightarrow w\mbox{ a.e. in }Q,\mbox{ as }\varepsilon \rightarrow 0,
\end{equation*}%
which implies that
\begin{equation*}
y_{\varepsilon }w_{\varepsilon }\rightarrow yw\mbox{ a.e. in }Q,\mbox{ as }%
\varepsilon \rightarrow 0.
\end{equation*}%
The functions $j$ and $j^{\ast }$ are continuous and so
\begin{equation*}
j(t,y_{\varepsilon }(t,x))\rightarrow j(t,y(t,x)),\mbox{ }j^{\ast
}(t,w_{\varepsilon }(t,x))\rightarrow j^{\ast }(t,w(t,x)),\mbox{ a.e. on }Q,%
\mbox{ as }\varepsilon \rightarrow 0.
\end{equation*}%
Now, by (\ref{68}) and (\ref{70-5}) we have
\begin{eqnarray*}
&&\int_{Q}\left( j(t,y_{\varepsilon }(t,x))+j^{\ast }(t,w_{\varepsilon
}(t,x))-y_{\varepsilon }w_{\varepsilon }\right) dxdt \\
&\leq &\int_{Q}\left( j(t,y(t,x))+j^{\ast }(t,w(t,x))-y_{\varepsilon
}w_{\varepsilon }\right) dxdt
\end{eqnarray*}%
and we apply the Fatou lemma because $j(t,y_{\varepsilon })+j^{\ast
}(t,w_{\varepsilon })-y_{\varepsilon }w_{\varepsilon }\geq 0.$ We get, using
(\ref{68}) and (\ref{70-5}) that
\begin{eqnarray*}
&&\int_{Q}\left( j(t,y(t,x))+j^{\ast }(t,w(t,x))-yw\right) dxdt \\
&\leq &\liminf\limits_{\varepsilon \rightarrow 0}\int_{Q}\left(
j(t,y_{\varepsilon }(t,x))+j^{\ast }(t,w_{\varepsilon }(t,x))-y_{\varepsilon
}w_{\varepsilon }\right) dxdt \\
&\leq &\limsup\limits_{\varepsilon \rightarrow 0}\int_{Q}\left(
j(t,x,y(t,x))+j^{\ast }(t,x,w(t,x))-y_{\varepsilon }w_{\varepsilon }\right)
dxdt \\
&\leq &\int_{Q}\left( j(t,y(t,x))+j^{\ast }(t,w(t,x))\right)
dxdt-\lim\limits_{\varepsilon \rightarrow 0}\int_{Q}y_{\varepsilon
}w_{\varepsilon }dxdt,
\end{eqnarray*}%
whence, by using (\ref{66}), we see that
\begin{equation}
-\int_{Q}ywdxdt\leq \frac{1}{2}\left\Vert y(T)\right\Vert _{V^{\prime }}^{2}-%
\frac{1}{2}\left\Vert y_{0}\right\Vert _{V^{\prime
}}^{2}-\int_{Q}yA_{0,\infty }^{-1}fdxdt.  \label{72}
\end{equation}%
We continue the proof by relying on the same arguments, starting this time
with Fatou's lemma applied for the positive function $j(t,x,-y_{\varepsilon
})+j^{\ast }(t,x,w_{\varepsilon })+y_{\varepsilon }w_{\varepsilon }$. By
similar computations we get
\begin{equation*}
-\int_{Q}ywdxdt\geq \frac{1}{2}\left\Vert y(T)\right\Vert _{V^{\prime }}^{2}-%
\frac{1}{2}\left\Vert y_{0}\right\Vert _{V^{\prime
}}^{2}-\int_{Q}yA_{0,\infty }^{-1}fdxdt,
\end{equation*}%
which together with (\ref{72}) imply (\ref{60}). \hfill $\square $

\medskip

Next we recall one of the main results given in \cite{gm-jota-var-1} in a
more general case, but particularized here to the space $L^{2}(Q).$

Let us consider the problem
\begin{eqnarray}
\frac{\partial y}{\partial t}-\Delta \widetilde{\beta }(t,x,y) &\ni &f\mbox{
\ \ \ \ \ \ \ \ \ \ \ \ \ \ in }Q,  \notag \\
-\frac{\partial \widetilde{\beta }(t,x,y)}{\partial \nu } &=&\alpha
\widetilde{\beta }(t,x,y)\mbox{ \ \ on }\Sigma ,  \label{si-1-0} \\
y(0,x) &=&y_{0}\mbox{ \ \ \ \ \ \ \ \ \ \ \ \ \ in }\Omega ,  \notag
\end{eqnarray}%
where $\widetilde{\beta }(t,x,r)=\partial \varphi (t,x,r)$ a.e. on $Q,$ for
all $r\in \mathbb{R},$ and $\varphi :\mathbb{R}\rightarrow \mathbb{R}$ is a
proper, convex, l.s.c. function satisfying $(h_{1}),(h_{2}),$ and the growth
condition
\begin{equation}
C_{1}\left\vert r\right\vert ^{2}+C_{1}^{0}\leq \varphi (t,x,r)\leq
C_{2}\left\vert r\right\vert ^{2}+C_{2}^{0},\mbox{ for all }r\in \mathbb{R},%
\mbox{ a.e. }t\in (0,T)  \label{si-1-1}
\end{equation}%
in addition. ($C_{i},$ $C_{i}^{0}$ are constants$,$ $i=1,2,$ and $C_{1}>0.)$

We consider the minimization problem
\begin{equation}
\mbox{Minimize }J_{0}(y,w)=\int_{Q}\left( \varphi (t,x,y(t,x))+\varphi
^{\ast }(t,x,w(t,x))-w(t,x)y(t,x)\right) dxdt  \tag{$P_{0}$}
\end{equation}%
for all $(y,w)\in U_{0},$ where
\begin{multline*}
U_{0}=\{(y,w);\mbox{ }y\in L^{2}(Q)\cap W^{1,2}([0,T];X_{2}^{\prime }),\mbox{
}y(T)\in V^{\prime },\mbox{ }w\in L^{2}(Q), \\
\varphi (\cdot ,\cdot ,y)\in L^{1}(Q),\mbox{ }\varphi ^{\ast }(\cdot ,\cdot
,w)\in L^{1}(Q), \\
(y,w)\mbox{ verifies (\ref{76}) below}\},
\end{multline*}%
\begin{eqnarray}
\frac{dy}{dt}(t)+\widetilde{A_{2}}w(t) &=&f(t)\mbox{ a.e. }t\in (0,T),
\label{76} \\
y(0) &=&y_{0}.  \notag
\end{eqnarray}%
We recall the notations $X_{2},$ $X_{2}^{\prime },$ $\widetilde{A_{2}}$
given in Section 2.1.

In \cite{gm-jota-var-1} it has been proved that $(P_{0})$ has at least a
solution and it has been established the equivalence between (\ref{si-1-0})
and $(P_{0}),$ resumed below (see Theorem 3.2 in \cite{gm-jota-var-1}).

\medskip

\noindent \textbf{Theorem 4.2.} \textit{Let }$y_{0}\in V^{\prime },$ $f\in
L^{\infty }(Q),$ \textit{and} \textit{let the pair }$(y,w)\in U_{0}$\textit{%
\ be a solution to }$(P_{0}).$\textit{\ Then,}%
\begin{equation*}
w(t,x)\in \widetilde{\beta }(t,y(t,x))\mbox{ \textit{a.e.} }(t,x)\in Q
\end{equation*}%
\textit{and }$(y,w)$\textit{\ is the unique weak solution to }(\ref{si-1-0}%
). \textit{Moreover,}%
\begin{eqnarray}
&&-\int_{0}^{t}\int_{\Omega }ywdxd\tau  \label{300} \\
&=&\frac{1}{2}\left\{ \left\Vert y(t)\right\Vert _{V^{\prime
}}^{2}-\left\Vert y_{0}\right\Vert _{V^{\prime }}^{2}\right\}
-\int_{0}^{t}\int_{\Omega }yA_{0,\infty }^{-1}fdxd\tau ,\mbox{ }for\mbox{ }%
all\mbox{ }t\in \lbrack 0,T].  \notag
\end{eqnarray}

\medskip

Of course, the result remains true when $\varphi $ does not depend on $x$.

Now, we can pass to the main result of this section which shows that a null
minimizer in $(P)$ provides a unique weak solution to (\ref{si-1}).

\medskip

\noindent \textbf{Theorem 4.3. }\textit{Under the assumptions }$%
(h_{1})-(h_{2}),$ (\ref{si-beta-j})-(\ref{si-9-2-0}), (\ref{si-9-7-01})-(\ref%
{si-9-7-02}), (\ref{si-9-5}) \textit{problem }$(P)$ \textit{has a solution }$%
(y^{\ast },w^{\ast })$\textit{\ such that }$y^{\ast }\in L^{\infty
}(0,T;V^{\prime }).$ \textit{Then, this solution is a null minimizer in (}$%
P) $\textit{\ }%
\begin{equation}
J(y^{\ast },w^{\ast })=\inf_{(y,w)\in U}J(y,w)=0  \label{Jinf}
\end{equation}%
\textit{and it turns out that it is the unique weak solution to }(\ref{si-1})%
\textit{. }

\medskip

\noindent \textbf{Proof. }Let us introduce the approximating problem
\begin{eqnarray}
\frac{\partial y}{\partial t}-\Delta \beta _{\lambda }(t,y) &=&f\mbox{ \ \ \
\ \ \ \ \ \ \ \ \ \ in }Q,  \notag \\
-\frac{\partial \beta _{\lambda }(t,y)}{\partial \nu } &=&\alpha \beta
_{\lambda }(t,y)\mbox{ \ \ on }\Sigma ,  \label{73} \\
y(0,x) &=&y_{0}\mbox{ \ \ \ \ \ \ \ \ \ \ \ \ in }\Omega ,  \notag
\end{eqnarray}%
where $\beta _{\lambda }$ is the Yosida approximation of $\beta .$

\noindent Let $\sigma $ be positive and consider the approximating problem
indexed upon $\sigma ,$%
\begin{eqnarray}
\frac{\partial y}{\partial t}-\Delta (\beta _{\lambda }(t,y)+\sigma y) &=&f%
\mbox{ \ \ \ \ \ \ \ \ \ \ \ \ \ \ \ \ \ \ \ in }Q,  \notag \\
-\frac{\partial (\beta _{\lambda }(t,y)+\sigma y)}{\partial \nu } &=&\alpha
(\beta _{\lambda }(t,y)+y)\mbox{ \ on }\Sigma ,  \label{74} \\
y(0,x) &=&y_{0}\mbox{ \ \ \ \ \ \ \ \ \ \ \ \ \ \ \ \ \ \ in }\Omega .
\notag
\end{eqnarray}%
The potential of $\beta _{\lambda }(t,r)+\sigma r$ is
\begin{equation}
j_{\lambda ,\sigma }(t,r)=j_{\lambda }(t,r)+\frac{\sigma }{2}r^{2},
\label{74-0}
\end{equation}%
where $j_{\lambda }$ is the Moreau regularization of $j.$ By a simple
computation using (\ref{70}), (\ref{si-9-7-01}), (\ref{si-9-2-0}) we get that%
\begin{equation}
\frac{\sigma }{2}\left\vert r\right\vert ^{2}+k_{1}r+k_{2}-2\lambda
k_{1}^{2}\leq j_{\lambda ,\sigma }(t,r)\leq j(t,0)+\left\vert r\right\vert
^{2}\left( \frac{1}{2\lambda }+\sigma ^{2}\right) .  \label{74-0-1}
\end{equation}%
Hence $j_{\lambda ,\sigma }$ satisfies (\ref{si-1-1}) and we rely on Theorem
4.2 with $\varphi (t,r)=j_{\lambda ,\sigma }(t,r)$ and $\widetilde{\beta }%
(t,r)=\beta _{\lambda }(t,r)+\sigma r$ to get that (\ref{74}) has a unique
weak solution $(y_{\lambda ,\sigma },w_{\lambda ,\sigma })\in U_{0},$
\begin{eqnarray*}
y_{\lambda ,\sigma } &\in &L^{2}(Q)\cap W^{1,2}([0,T];X_{2}^{\prime }),\mbox{
}y_{\lambda ,\sigma }(T)\in V^{\prime }, \\
w_{\lambda ,\sigma } &=&\beta _{\lambda }(t,y_{\lambda ,\sigma })+\sigma
y_{\lambda ,\sigma }\in L^{2}(Q).
\end{eqnarray*}%
This solution is the null minimizer in $(P_{0}),$ i.e.,
\begin{equation}
J_{0}(y_{\lambda ,\sigma },w_{\lambda ,\sigma })=\int_{Q}\left( j_{\lambda
,\sigma }(t,y_{\lambda ,\sigma }(t,x))+j_{\lambda ,\sigma }^{\ast
}(t,w_{\lambda ,\sigma }(t,x))-y_{\lambda ,\sigma }w_{\lambda ,\sigma
}\right) dxdt=0,  \label{77}
\end{equation}%
and satisfies (\ref{76}), namely%
\begin{eqnarray}
\frac{dy_{\lambda ,\sigma }}{dt}(t)+\widetilde{A_{2}}w_{\lambda ,\sigma }(t)
&=&f(t)\mbox{ a.e. }t\in (0,T),  \label{76-0} \\
y(0) &=&y_{0}.  \notag
\end{eqnarray}%
Moreover, we have by (\ref{300}) that
\begin{eqnarray}
&&-\int_{0}^{t}\int_{\Omega }y_{\lambda ,\sigma }w_{\lambda ,\sigma }dxd\tau
\label{si-9-1} \\
&=&\frac{1}{2}\left\{ \left\Vert y_{\lambda ,\sigma }(t)\right\Vert
_{V^{\prime }}^{2}-\left\Vert y_{0}\right\Vert _{V^{\prime }}^{2}\right\}
-\int_{0}^{t}\int_{\Omega }y_{\lambda ,\sigma }A_{0,\infty }^{-1}fdxd\tau ,%
\mbox{ }  \notag
\end{eqnarray}%
for all $t\in \lbrack 0,T].$ Taking into account (\ref{si-9-1}) and (\ref{77}%
) we still can write
\begin{eqnarray}
&&\int_{Q}(j_{\lambda ,\sigma }(t,y_{\lambda ,\sigma }(t,x))+j_{\lambda
,\sigma }^{\ast }(t,w_{\lambda ,\sigma }(t,x)))dxdt  \label{78} \\
&&+\frac{1}{2}\left\{ \left\Vert y_{\lambda ,\sigma }(T)\right\Vert
_{V^{\prime }}^{2}-\left\Vert y_{0}\right\Vert _{V^{\prime }}^{2}\right\}
=\int_{0}^{T}\int_{\Omega }y_{\lambda ,\sigma }A_{0,\infty }^{-1}fdxdt.
\notag
\end{eqnarray}%
We note that
\begin{equation}
\frac{j_{\lambda ,\sigma }^{\ast }(t,\omega )}{\left\vert \omega \right\vert
}\rightarrow \infty \mbox{ as }\left\vert \omega \right\vert \rightarrow
\infty ,  \label{78-0}
\end{equation}%
uniformly in $\lambda $ and $\sigma .$ This happens due to (\ref{si-9-7})
because by setting
\begin{equation*}
\eta _{\lambda ,\,\sigma }=\partial j_{\lambda ,\sigma }(t,r),\mbox{ }\eta
_{\lambda ,\sigma }=\beta _{\lambda }(t,r)+\sigma r=\beta (t,(1+\lambda
\beta (t,\cdot ))^{-1}r)+\sigma r,
\end{equation*}%
then $\eta _{\lambda ,\sigma }$ is bounded on bounded subsets $\left\vert
r\right\vert \leq M,$ uniformly in $\lambda $ and $\sigma ,$ for $\lambda $
and $\sigma $ small (smaller than 1, e.g.).

We also note that
\begin{eqnarray*}
&&\int_{Q}j_{\lambda ,\sigma }(t,y_{\lambda ,\sigma }(t,x))dxdt\geq
\int_{Q}j_{\lambda }(t,y_{\lambda ,\sigma }(t,x))dxdt \\
&=&\int_{Q}j(t,(1+\lambda \beta (t,\cdot ))^{-1}y_{\lambda ,\sigma
})dxdt+\int_{Q}\frac{1}{2\lambda }\left\vert y_{\lambda ,\sigma }-(1+\lambda
\beta (t,\cdot ))^{-1}y_{\lambda ,\sigma }\right\vert ^{2}dxdt
\end{eqnarray*}%
and
\begin{eqnarray*}
&&\int_{0}^{T}\int_{\Omega }y_{\lambda ,\sigma }A_{0,\infty
}^{-1}fdxdt=\int_{Q}\left( y_{\lambda ,\sigma }-(1+\lambda \beta (t,\cdot
))^{-1}y_{\lambda ,\sigma }\right) A_{0,\infty }^{-1}fdxdt \\
&&+\int_{Q}(1+\lambda \beta (t,\cdot ))^{-1}y_{\lambda ,\sigma }A_{0,\infty
}^{-1}fdxdt \\
&\leq &\int_{Q}\frac{1}{2\lambda }\left\vert y_{\lambda ,\sigma }-(1+\lambda
\beta (t,\cdot ))^{-1}y_{\lambda ,\sigma }\right\vert ^{2}dxdt+2\lambda
\int_{Q}(A_{0,\infty }^{-1}f)^{2}dxdt \\
&&+\int_{Q}(1+\lambda \beta (t,\cdot ))^{-1}y_{\lambda ,\sigma }A_{0,\infty
}^{-1}fdxdt.
\end{eqnarray*}%
Plugging these in (\ref{78}) we get after some algebra that
\begin{eqnarray}
&&\int_{Q}j(t,(1+\lambda \beta (t,\cdot ))^{-1}y_{\lambda ,\sigma
})dxdt+\int_{Q}j_{\lambda ,\sigma }^{\ast }(t,w_{\lambda ,\sigma }(t,x)))dxdt
\label{78-1} \\
&&+\frac{1}{2}\left\{ \left\Vert y_{\lambda ,\sigma }(T)\right\Vert
_{V^{\prime }}^{2}-\left\Vert y_{0}\right\Vert _{V^{\prime }}^{2}\right\}
\notag \\
&\leq &\int_{Q}(1+\lambda \beta (t,\cdot ))^{-1}y_{\lambda ,\sigma
}A_{0,\infty }^{-1}fdxdt+2\lambda \int_{Q}(A_{0,\infty }^{-1}f)^{2}dxdt.
\notag
\end{eqnarray}%
Further we set
\begin{equation*}
(1+\lambda \beta (t,\cdot ))^{-1}y_{\lambda ,\sigma }=z_{\lambda ,\sigma }
\end{equation*}%
and argue as in Theorem 3.2 to deduce by the Dunford-Pettis theorem that $%
(z_{\lambda ,\sigma })_{\sigma }$ and $(w_{\lambda ,\sigma })_{\sigma }$ are
weakly compact in $L^{1}(Q).$ Recalling (\ref{bound}) we also get%
\begin{equation}
\left\Vert y_{\lambda ,\sigma }(T)\right\Vert _{V^{\prime }}\leq C
\label{79}
\end{equation}%
independently on $\sigma $ and $\lambda $.

Taking into account that $w_{\lambda ,\sigma }=\beta _{\lambda }(y_{\lambda
,\sigma })+\sigma y_{\lambda ,\sigma },$ equation (\ref{si-9-1}) yields
\begin{eqnarray*}
&&\frac{1}{2}\left\Vert y_{\lambda ,\sigma }(t)\right\Vert _{V^{\prime
}}^{2}+\int_{0}^{t}\int_{\Omega }(\beta _{\lambda }(\tau ,y_{\lambda ,\sigma
})y_{\lambda ,\sigma }+\sigma y_{\lambda ,\sigma }^{2})dxd\tau \\
&=&\frac{1}{2}\left\Vert y_{0}\right\Vert _{V^{\prime
}}^{2}+\int_{0}^{t}\left\langle y_{\lambda ,\sigma }(\tau ),A_{0,\infty
}^{-1}f(\tau )\right\rangle _{V^{\prime },V}d\tau \\
&=&\frac{1}{2}\left\Vert y_{0}\right\Vert _{V^{\prime }}^{2}+\frac{1}{2}%
\int_{0}^{t}\left\Vert A_{0,\infty }^{-1}f(\tau )\right\Vert _{V}^{2}d\tau +%
\frac{1}{2}\int_{0}^{t}\left\Vert y_{\lambda ,\sigma }(\tau )\right\Vert
_{V^{\prime }}^{2}d\tau ,\mbox{ }
\end{eqnarray*}%
for all $t\in \lbrack 0,T].$

Taking into account that $\beta _{\lambda }(t,r)r\geq 0,$ for all $r\in
\mathbb{R},$ and in virtue of the Gronwall lemma, we deduce that
\begin{equation}
\left\Vert y_{\lambda ,\sigma }\right\Vert _{L^{\infty }(0,T;V^{\prime
})}\leq C,  \label{79-0}
\end{equation}%
\begin{equation*}
\sqrt{\sigma }\left\Vert y_{\lambda ,\sigma }\right\Vert _{L^{2}(Q)}\leq C,
\end{equation*}%
and%
\begin{equation}
\int_{Q}j(t,z_{\lambda ,\sigma }(t,x))dxdt\leq C,\mbox{ }\int_{Q}j_{\lambda
,\sigma }^{\ast }(t,w_{\lambda ,\sigma }(t,x))dxdt\leq C  \label{79-1}
\end{equation}%
independently on $\sigma $ and $\lambda .$ (For getting (\ref{79-1}) we
recall the arguments leading to (\ref{101-1}), (\ref{101-2})).

Then, (\ref{78}) and relation (\ref{si-9-7-01}) for $j_{\lambda ,\sigma
}^{\ast }$ imply that
\begin{equation}
\int_{Q}j_{\lambda ,\sigma }(t,y_{\lambda ,\sigma }(t,x))dxdt\leq C
\label{79-2}
\end{equation}%
independently on $\sigma $ and $\lambda .$ Following the proof of Theorem
3.2 we deduce that
\begin{eqnarray*}
z_{\lambda ,\sigma } &\rightarrow &z_{\lambda }\mbox{ weakly in }L^{1}(Q),%
\mbox{ as }\sigma \rightarrow 0, \\
w_{\lambda ,\sigma } &\rightarrow &w_{\lambda }\mbox{ weakly in }L^{1}(Q),%
\mbox{ as }\sigma \rightarrow 0, \\
\sqrt{\sigma }y_{\lambda ,\sigma } &\rightarrow &\zeta _{\lambda }\mbox{
weakly in }L^{2}(Q),\mbox{ as }\sigma \rightarrow 0, \\
y_{\lambda ,\sigma } &\rightarrow &y_{\lambda }\mbox{ weak-star in }%
L^{\infty }(0,T;V^{\prime }),\mbox{ as }\sigma \rightarrow 0, \\
y_{\lambda ,\sigma }(T) &\rightarrow &\xi \mbox{ weakly in }V^{\prime },%
\mbox{ as }\sigma \rightarrow 0, \\
Aw_{\lambda } &\rightarrow &Aw_{\lambda }\mbox{ weakly in }%
L^{1}(0,T;X^{\prime }),\mbox{ as }\sigma \rightarrow 0, \\
\frac{dy_{\lambda ,\sigma }}{dt} &\rightarrow &\frac{dy_{\lambda }}{dt}\mbox{
weakly in }L^{1}(0,T;X^{\prime }),\mbox{ as }\sigma \rightarrow 0.
\end{eqnarray*}%
By (\ref{79-2}) and (\ref{70-0}) we have
\begin{equation*}
\int_{Q}\frac{1}{2\lambda }\left\vert y_{\lambda ,\sigma }-(1+\lambda \beta
(t,\cdot ))^{-1}y_{\lambda ,\sigma }\right\vert ^{2}dxdt\leq
\int_{Q}j_{\lambda ,\sigma }(t,y_{\lambda ,\sigma }(t,x))dxdt\leq C
\end{equation*}%
whence, denoting $\chi _{\lambda ,\sigma }=\left( y_{\lambda ,\sigma
}-z_{\lambda ,\sigma }\right) /\sqrt{2\lambda }$ we see that $(\chi
_{\lambda ,\sigma })_{\sigma }$ is bounded in $L^{2}(Q)$ and $\chi _{\lambda
,\sigma }\rightarrow \chi _{\lambda }$ weakly in $L^{2}(Q),$ as $\sigma
\rightarrow 0,$ on a subsequence. Then
\begin{equation*}
y_{\lambda ,\sigma }-z_{\lambda ,\sigma }\rightarrow \sqrt{2\lambda }\chi
_{\lambda }\mbox{ weakly in }L^{1}(Q),\mbox{ as }\sigma \rightarrow 0,
\end{equation*}%
where $\left\Vert \chi _{\lambda }\right\Vert _{L^{1}(Q)}\leq C.$ Since $%
z_{\lambda ,\sigma }\rightarrow z_{\lambda }$ weakly in $L^{1}(Q),$ it
follows that $(y_{\lambda ,\sigma })_{\sigma }$ is bounded in $L^{1}(Q),$ so
it converges weakly and by the limit uniqueness we have
\begin{equation*}
y_{\lambda ,\sigma }\rightarrow y_{\lambda }\mbox{ weakly in }L^{1}(Q),\mbox{
as }\sigma \rightarrow 0.
\end{equation*}%
We also have%
\begin{equation}
y_{\lambda }=z_{\lambda }+\sqrt{2\lambda }\zeta _{\lambda }\mbox{ a.e. on }Q.
\label{79-3}
\end{equation}%
By Arzel\`{a}-Ascoli theorem (since $V^{\prime }$ is compact in $X^{\prime }$
because $X$ is compact in $V)$ it follows that
\begin{equation*}
y_{\lambda ,\sigma }(t)\rightarrow y_{\lambda }(t)\mbox{ in }X^{\prime },%
\mbox{ uniformly in }t\in \lbrack 0,T],\mbox{ as }\sigma \rightarrow 0,
\end{equation*}%
so $\xi =y_{\lambda }(T)$ and $y_{\lambda }(0)=y_{0}.$

Passing to the limit in (\ref{76-0}) we get that $(y_{\lambda },w_{\lambda
}) $ satisfies
\begin{eqnarray}
\frac{dy_{\lambda }}{dt}(t)+Aw_{\lambda }(t) &=&f(t)\mbox{ a.e. }t\in (0,T),
\label{76-00} \\
y(0) &=&y_{0}.  \notag
\end{eqnarray}

Passing to the limit in (\ref{78-1}) as $\sigma \rightarrow 0$, using the
weak lower semicontinuity property we get
\begin{eqnarray}
\int_{Q}(j(t,z_{\lambda }(t,x))+j_{\lambda }^{\ast }(t,w_{\lambda
}(t,x)))dxdt &&  \label{82} \\
+\frac{1}{2}\left\{ \left\Vert y_{\lambda }(T)\right\Vert _{V^{\prime
}}^{2}-\left\Vert y_{0}\right\Vert _{V^{\prime }}^{2}\right\}
-\int_{Q}y_{\lambda }A_{0,\infty }^{-1}fdxdt-2\lambda \int_{Q}(A_{0,\infty
}^{-1}f)^{2}dxdt &\leq &0.  \notag
\end{eqnarray}

We repeat again the arguments developed in Theorem 3.2 and deduce by the
Dunford-Pettis theorem that $(z_{\lambda })_{\lambda }$ and $(w_{\lambda
})_{\lambda }$ are weakly compact in $L^{1}(Q).$ It still follows that
\begin{equation}
\left\Vert z_{\lambda }(T)\right\Vert _{V^{\prime }}\leq C,\mbox{ }%
\int_{Q}j(t,z_{\lambda }(t,x))dxdt\leq C,\mbox{ }\int_{Q}j_{\lambda }^{\ast
}(t,w_{\lambda }(t,x))dxdt\leq C  \label{82-1}
\end{equation}%
independently on $\lambda $ (recall (\ref{bound}), (\ref{101-1}), (\ref%
{101-2}))$.$ Passing to the limit in (\ref{79-0}) as $\sigma \rightarrow 0$
we get
\begin{equation*}
\left\Vert y_{\lambda }\right\Vert _{L^{\infty }(0,T;V^{\prime })}\leq C
\end{equation*}%
where $C$ are several constants independent on $\lambda .$

Then, proceeding along with the proof of Theorem 3.2 we obtain from (\ref{82}%
), by selecting subsequences, that
\begin{eqnarray*}
z_{\lambda } &\rightarrow &z^{\ast }\mbox{ weakly in }L^{1}(Q),\mbox{ as }%
\lambda \rightarrow 0, \\
w_{\lambda } &\rightarrow &w^{\ast }\mbox{ weakly in }L^{1}(Q),\mbox{ as }%
\sigma \rightarrow 0, \\
y_{\lambda } &\rightarrow &y^{\ast }\mbox{ weak-star in }L^{\infty
}(0,T;V^{\prime }),\mbox{ as }\sigma \rightarrow 0, \\
y_{\lambda }(T) &\rightarrow &y^{\ast }(T)\mbox{ weakly in }V^{\prime },%
\mbox{ as }\sigma \rightarrow 0, \\
Aw_{\lambda } &\rightarrow &Aw^{\ast }\mbox{ weakly in }L^{1}(0,T;X^{\prime
}),\mbox{ as }\sigma \rightarrow 0, \\
\frac{dy_{\lambda }}{dt} &\rightarrow &\frac{dy^{\ast }}{dt}\mbox{ weakly in
}L^{1}(0,T;X^{\prime }),\mbox{ as }\sigma \rightarrow 0.
\end{eqnarray*}%
By (\ref{79-3}) we get that $z^{\ast }=y^{\ast }$ a.e. on $Q$ and by (\ref%
{82-1}) we obtain%
\begin{equation}
\left\Vert y^{\ast }(T)\right\Vert _{V^{\prime }}\leq C,\mbox{ }%
\int_{Q}j(t,y^{\ast }(t,x))dxdt\leq C,\mbox{ }\int_{Q}j^{\ast }(t,w^{\ast
}(t,x))dxdt\leq C.  \label{82-3}
\end{equation}%
The first inequality is obvious. For the second (if $j(t,r)\geq 0)$ Fatou's
lemma yields
\begin{equation}
\int_{Q}j(t,y^{\ast }(t,x))dxdt=\int_{Q}\liminf\limits_{\lambda \rightarrow
0}j(t,z_{\lambda }(t,x))dxdt\leq \liminf\limits_{\lambda \rightarrow
0}\int_{Q}j(t,z_{\lambda }(t,x))dxdt\leq C.  \label{82-2}
\end{equation}%
If $j$ is not positive, we use again (\ref{si-9-7-01}) and denoting $%
\widetilde{j}(t,r)=j(t,r)-k_{1}r-k_{2}\geq 0$ we write
\begin{equation*}
\int_{Q}\widetilde{j}(t,y^{\ast }(t,x))dxdt\leq \liminf\limits_{\lambda
\rightarrow 0}\int_{Q}\widetilde{j}(t,z_{\lambda }(t,x))dxdt,
\end{equation*}%
whence we get
\begin{eqnarray*}
&&\int_{Q}j(t,y^{\ast }(t,x))dxdt-\int_{Q}(k_{1}y^{\ast }+k_{2})dxdt \\
&\leq &\liminf\limits_{\lambda \rightarrow 0}\int_{Q}j(t,z_{\lambda
}(t,x))dxdt-\int_{Q}(k_{1}y^{\ast }+k_{2})dxdt,
\end{eqnarray*}%
i.e., (\ref{82-2}). In what concerns the third inequality in (\ref{82-3}),
we can write by (\ref{82-1})
\begin{eqnarray}
&&\int_{Q}\frac{1}{2\lambda }\left\vert w_{\lambda }-(1+\lambda \beta
^{-1}(t,\cdot ))^{-1}w_{\lambda }\right\vert ^{2}dxdt+j^{\ast }(t,(1+\beta
^{-1}(t,\cdot ))^{-1}w_{\lambda })  \label{82-4} \\
&\leq &\int_{Q}j_{\lambda }^{\ast }(t,w_{\lambda }(t,x))dxdt\leq C.  \notag
\end{eqnarray}%
Recalling that $j^{\ast }(t,x,\omega )\geq k_{3}(t,x)\omega +k_{4}(t,x)$ we
get that%
\begin{equation*}
k_{3}\int_{Q}(1+\beta ^{-1}(t,\cdot ))^{-1}w_{\lambda }dxdt\leq C+k_{4}\mbox{%
meas}(Q),
\end{equation*}%
hence $\left( (1+\beta ^{-1}(t,\cdot ))^{-1}w_{\lambda }\right) _{\lambda }$
is bounded in $L^{1}(Q).$ Then%
\begin{equation*}
\int_{Q}\frac{1}{2\lambda }\left\vert w_{\lambda }-(1+\lambda \beta
^{-1}(t,\cdot ))^{-1}w_{\lambda }\right\vert ^{2}dxdt\leq
C+\int_{Q}(k_{3}(1+\beta ^{-1}(t,\cdot ))^{-1}w_{\lambda }+k_{4})dxdt\leq
C_{1}.
\end{equation*}%
It follows that
\begin{equation*}
(1+\lambda \beta ^{-1}(t,\cdot ))^{-1}w_{\lambda }\rightarrow w^{\ast }\mbox{
weakly in }L^{1}(Q),\mbox{ as }\lambda \rightarrow 0.
\end{equation*}%
Then, we passing to the limit in (\ref{82-4}) as $\lambda \rightarrow 0$ (if
$j^{\ast }$ is nonnegative). Otherwise we use again (\ref{si-9-7-01}) for $%
j^{\ast }.$

Passing to the limit in (\ref{76-00}) and (\ref{82}) as $\lambda \rightarrow
0$ we get%
\begin{eqnarray}
\frac{dy^{\ast }}{dt}(t)+Aw^{\ast }(t) &=&f(t)\mbox{ a.e. }t\in (0,T),
\label{76-000} \\
y(0) &=&y_{0},  \notag
\end{eqnarray}%
and, again by the weak lower semicontinuity,
\begin{eqnarray}
&&\int_{Q}(j(t,y^{\ast }(t,x))+j^{\ast }(t,w^{\ast }(t,x)))dxdt  \label{83}
\\
&&+\frac{1}{2}\left\{ \left\Vert y^{\ast }(T)\right\Vert _{V^{\prime
}}^{2}-\left\Vert y_{0}\right\Vert _{V^{\prime }}^{2}\right\}
-\int_{Q}y^{\ast }A_{0,\infty }^{-1}fdxd\tau \leq 0.  \notag
\end{eqnarray}%
We have got that $(y^{\ast },w^{\ast })\in U,$ $y^{\ast }\in L^{\infty
}(0,T;V^{\prime })$ and so by Lemma 4.1 it follows that $y^{\ast }w^{\ast
}\in L^{1}(Q).$ Replacing the sum of the last two terms on the right-hand
side in (\ref{83}) by (\ref{60}) we get
\begin{equation*}
\int_{Q}(j(t,y^{\ast }(t,x))+j^{\ast }(t,w^{\ast }(t,x))-y^{\ast
}(t,x)w^{\ast }(t,x))dxdt\leq 0.
\end{equation*}%
Recalling (\ref{si-4-1}) we obtain
\begin{equation}
\int_{Q}(j(t,y^{\ast }(t,x))+j^{\ast }(t,w^{\ast }(t,x))-y^{\ast
}(t,x)w^{\ast }(t,x))dxdt=0  \label{84}
\end{equation}%
which eventually implies that
\begin{equation*}
j(t,y^{\ast }(t,x))+j^{\ast }(t,w^{\ast }(t,x))-y^{\ast }(t,x)w^{\ast
}(t,x)=0\mbox{ a.e. on }Q.
\end{equation*}%
Therefore, we conclude that $w^{\ast }(t,x)\in \beta (t,y^{\ast }(t,x))$
a.e. on $Q,$ by the Legendre-Fenchel relations.

On the other hand, due again to (\ref{60}) in Lemma 4.1, relation (\ref{84})
means in fact that $(y^{\ast },w^{\ast })$ realizes the minimum in $(P),$ as
claimed in (\ref{Jinf}).

The uniqueness follows directly by (\ref{si-1}) using the monotony of $\beta
.$

Indeed, let $(y$,$\eta )$ and $(\widetilde{y},\widetilde{\eta })$ be two
solutions to (\ref{si-1}) corresponding to the same data, where $\eta
(t,x)\in \beta (t,x,y(t,x)),$ $\widetilde{\eta }(t,x)\in \beta (t,x,%
\widetilde{y}(t,x))$ a.e. on $Q,$ $(y,\eta )$ and $(\widetilde{y},\widetilde{%
\eta })$ belong to $U$ and $y,\widetilde{y}\in L^{\infty }(0,T;V^{\prime }).$
We write the equations satisfied by their difference%
\begin{eqnarray*}
\frac{d(y-\widetilde{y})}{dt}(t)+A(\eta -\widetilde{\eta })(t) &\ni &0\mbox{
a.e. }t\in (0,T), \\
(y-\widetilde{y})(0) &=&0,
\end{eqnarray*}%
multiply the equation by $A_{0,\infty }^{-1}(y-\widetilde{y})(t)$ and
integrate it over $(0,t)$ obtaining
\begin{equation*}
\frac{1}{2}\left\Vert (y-\widetilde{y})(t)\right\Vert _{V^{\prime
}}^{2}+\int_{0}^{t}\int_{\Omega }\left( \eta -\widetilde{\eta })(y-%
\widetilde{y}\right) dxdt=0.
\end{equation*}%
But $\beta (t,r)$ is maximal monotone, hence we get $\left\Vert y(t)-%
\widetilde{y}(t)\right\Vert _{V^{\prime }}^{2}\leq 0,$ whence $y(t)=%
\widetilde{y}(t)$ for all $t\in \lbrack 0,T].$ \hfill $\square $

\bigskip

We remark now that if $(y^{\ast },\eta ^{\ast })$ is the solution to (\ref%
{si-1}), with $\eta ^{\ast }(t,x)\in \beta (t,y^{\ast }(t,x))$ a.e. on $Q$
and $y^{\ast }\in L^{\infty }(0,T;V^{\prime }),$ then it is a unique
solution to $(P),$ because by Lemma 4.1, $J(y,w)\geq 0$ for any $(y,w)\in U$
and the minimum is realized at $(y^{\ast },\eta ^{\ast })$ since $J(y^{\ast
},\eta ^{\ast })=0.$ So, we conclude that (\ref{si-1}) is equivalent with
the minimization problem $(P)$.

\section{Conclusions}

This paper deals with the application of the Brezis-Ekeland principle for a
nonlinear diffusion equation with a monotonically increasing time depending
nonlinearity which provides a potential having a weak coercivity property.
The result states that the solution of the nonlinear equation can be
retrieved as the null minimizer of an appropriate minimization problem for a
convex functional involving the potential of the nonlinearity.

This approach is useful because it allows the existence proof in cases in
which, due to the generality of the nonlinearity, standard methods do not
apply. Also it can lead to a simpler numerical computation of the solution
to the equation by replacing its direct determination by the numeric
calculus of the minimum of a convex functional with a linear state equation.

With respect to the literature concerning existence results for (\ref{si-1}%
), Theorem 4.3 provides existence under very general conditions on the
nonlinear function $\beta .$ As regards the assumption (\ref{si-9-5}) it can
be equivalently expressed as
\begin{equation*}
\limsup_{\left\vert r\right\vert \rightarrow \infty }\frac{j(t,-r)}{j(t,r)}%
<\infty ,
\end{equation*}%
(see \cite{barbu-daprato}). Since in specific real problems the solution to (%
\ref{si-1}) is nonnegative, and so $\beta $ is defined on $[0,\infty ),$
this condition is achieved by extending in a convenient way the function $%
\beta $ on $(-\infty ,0).$ For instance, conditions (\ref{si-9-2})-(\ref%
{si-9-2-00})\ are satisfied for $\beta $ of the form
\begin{equation*}
\beta (t,x,r)=\mbox{sgn(}r)\log (\left\vert r\right\vert +a(t,x)),\mbox{ }%
a\geq a_{0}>0
\end{equation*}%
or
\begin{equation*}
\beta (t,x,r)=\mbox{sgn}(r)\exp (a(t,x)r^{2}),\mbox{ }a\geq a_{0}>0.
\end{equation*}

Concerning possible applications, we remark that problem (\ref{si-1}) can be
obtained by a change of variable in an equation of the form
\begin{equation*}
\frac{\partial (m(t,x)y)}{\partial t}-\Delta \beta _{0}(y)\ni f
\end{equation*}%
which is associated to various physical models as for example: fluid
diffusion in saturated-unsaturated deformable porous media with the porosity
$m$ time and space dependent (see appropriate problems in \cite{af-gm-12},
\cite{gm-cc-08}), or to absorption-desorption processes in saturated porous
media in which $m$ is the absorption-desorption rate of the fluid by the
solid. The Robin boundary condition arising in (\ref{si-1}) was chosen
because of its relevance in these physical models. Also, evolution equations
with nonautonomous operators can be associated to models in which the
boundary conditions are of time dependent nonhomogeneous Dirichlet type, or
nonlocal as in population dynamics (see \cite{cuiama-2}).

In all these problems the coefficient $m$ or the coefficients in the
boundary conditions may have a very low regularity which makes not possible
the approach of (\ref{si-1}) by the nonlinear semigroup method in the
time-dependent case given in \cite{Crandall-Pazy-71}.

\bigskip

\noindent \textbf{Acknowledgment. }This work has been carried out within a
grant of the Romanian National Authority for Scientific Research,
CNCS-UEFISCDI, project number PN-II-ID-PCE-2011-3-0027. The author
acknowledges the fellowship awarded by CIRM-Fondazione Bruno Kessler, Italy,
in November 2012.

\end{document}